\magnification=1200
\baselineskip =13pt
\overfullrule =0pt
\ifnum\pageno=1\nopagenumbers

\def\op{\operatorname}

\input amstex
\documentstyle{amsppt}
\nologo
\vsize=500pt

\topmatter

\title{Derived quot schemes}\endtitle
\author{Ionu\c t Ciocan-Fontanine and
Mikhail M. Kapranov}\endauthor

\address{Department of Mathematics,
Northwestern University, Evanston, IL 60208}\endaddress

\email ciocan\@math.nwu.edu,  kapranov\@math.nwu.edu \endemail

\endtopmatter

\document

\heading {Introduction}\endheading

\noindent {\bf (0.1)} A typical moduli problem in geometry
is to construct a ``space" $H$ parametrizing, up to isomorphism,
objects
of some given category $\Cal  Z$ (e.g., manifolds, vector bundles etc.).
This can be seen as a kind of a non-Abelian cohomology problem and
the construction usually consists of two steps of opposite nature, namely
applying a left exact functor (A) followed by a right
exact functor (B):

\vskip .2cm

\noindent (A) One finds a space $Z$ of ``cocycles" whose points parametrize
objects of $\Cal  Z$ equipped with some extra structure. Usually $Z$
is given inside a much simpler space $C$ of ``cochains" by
explicit equations, so forming $Z$ is an inverse limit-type
construction (hence left exact).

\vskip .2cm

\noindent (B) One factorizes $Z$ by the action of a group (or groupoid) $B$
by identifying isomorphic objects and sets $H = Z/B$. This is a direct
limit-type construction, hence right exact. 

\vskip .2cm

Part (B) leads to well known difficulties which in algebraic geometry
are resolved by using the language of stacks. This can be seen
as passing to the nonabelian left derived functor of (B). Indeed,
an algebraic stack is a nonlinear analog of a complex of vector
spaces situated in degrees $[-1,0]$ and, for example, the
tangent ``space" to a stack at a point is a complex of this nature. 

  The step (A) may or may not be as clearly noticeable
because points of $C$ have no meaning from the point of view of the category
$\Cal  Z$. It is also very important, nevertheless, because $C$ is usually
smooth while $Z$ and hence $H$ may be singular (even as a stack). 
For example, when $\Cal  Z$ consists of complex analytic vector
bundles, we can take $Z$ to consist of integrable $\bar\partial$-connections
on a given smooth bundle. Then $C$ consists of all 
 $\bar\partial$-connections, integrable or not, which form an
infinite-dimensional affine space, but do not, in general, define
holomorphic bundles.

\vskip .3cm

\noindent {\bf (0.2)} The derived deformation theory (DDT) program,
see \cite{Kon} \cite{Ka1} \cite{Hi}, is a program of research aimed at 
systematically
resolving the difficulties related to singularities of the
moduli spaces. It is convenient to formulate its most important
premises as follows:

\vskip .1cm

\noindent {\bf (a)} One should take the right derived functor in the step (A) 
as well, landing
in an appropriate ``right derived category of schemes" whose objects
(called dg-schemes)
are nonlinear analogs of cochain complexes situated in degrees
$[0,\infty)$ and whose tangent spaces are indeed complexes of this kind.

\vskip .1cm

\noindent {\bf (b)} The object $RZ$ obtained in this way, should be manifestly 
smooth
in an appropriate sense (so that the singular nature of $Z$
is the result of truncation).

\vskip .1cm

\noindent {\bf (c)} The correct moduli ``space" $LRH$ is the result of 
applying to $RZ$
the (stack-theoretic) left derived functor of (B). It should
lie in a larger derived category of ``dg-stacks" whose linear
objects are cochain complexes situated in degrees $[-1, +\infty)$.
The tangent space  to $LRH$ at a point corresponding
to an object $X\in\Cal  Z$  is a complex of this kind, and its
$i$th cohomology space is naturally identified to the
$(i+1)$st cohomology space of the sheaf of infinitesimal automorphisms of $X$, 
thus generalizing
the classical Kodaira-Spencer isomorphism to higher
cohomology. 

\vskip .1cm

\noindent {\bf (d)} All considerations in algebraic geometry
which involve deformation to a generic almost
complex structure can and should
be replaced by systematically working with the derived moduli space
$LRH$, its characteristic classes etc.  

\vskip .1cm

It is important not to confuse the putative dg-stacks of (c)
with algebraic $n$-stacks as developed by Simpson \cite{Si}:
the latter serve as nonlinear analogs of cochain complexes
situated in degrees $[-n, 0]$.

\vskip .3cm

\noindent {\bf (0.3)} In this paper we concentrate on
taking the derived functor of the step (A) in
the framework of algebraic geometry. 
Classically, almost all  constructions of moduli spaces
in this framework proceed via Hilbert schemes and their
generalizations, $Quot$ schemes, introduced by 
Grothendieck \cite{Gr}, see \cite{Kol} \cite{Vi}
for detailed exposition. 
 In many cases, the construction
goes simply by quotienting an appropriate part of the $Quot$ scheme
by an action of an algebraic group, thus giving an algebraic stack.
(Thus the scheme $Quot$ plays the role of $Z$ in (0.1)(A)).

The first step in
constructing
derived moduli spaces in algebraic geometry is, then, to
construct the derived version of $Quot$. This is done in
the present paper.
 
To recall the situation, let $\Bbb K$ be a field,
 let $X$ be a projective scheme over $\Bbb K$ and $\Cal  F$ be
a coherent sheaf on $X$. The scheme $Quot(\Cal  F)$  can be viewed as 
parametrizing coherent subsheaves $\Cal  K\subset\Cal  F$:
to every such $\Cal  K$ there corresponds a $\Bbb K$-point $[\Cal  K]
\in Quot(\Cal  F)$.

Assume that $\Bbb K$ has
characteristic 0.  For each Hilbert polynomial $h$
we construct a smooth dg-manifold (see \S 2 for background)
$RQuot_h(\Cal  F)$ 
 with the following properties:

\vskip .1cm

\noindent (0.3.1) The degree 0 truncation $\pi_0 RQuot_h(\Cal  F)$ is 
identified
with $Quot_h(\Cal  F)$.

\vskip .1cm

\noindent (0.3.2) If $[\Cal  K]$ is a $\Bbb K$-point of $Quot_h(\Cal  F)$
corresponding to a subsheaf $\Cal  K\subset \Cal  F$, then the
tangent space $T^\bullet_{[\Cal  K]}RQuot_h(\Cal  F)$ is a
$\Bbb Z_+$-graded cochain complex whose cohomology is given by:
$$H^i T^\bullet_{[\Cal  K]}RQuot_h(\Cal  F)\simeq \text{Ext}^i_{\Cal  
O_X}(\Cal  K,
\Cal  F/\Cal  K).$$

\vskip .1cm

Note that for the ordinary Quot scheme the tangent space is given by
taking $i=0$ in (0.3.2) (i.e., by $\text{Hom}_{\Cal  O_X}(\Cal  K, \Cal  
F/\Cal  K)$).
 It is perfectly possible for the dimension of this Hom to jump in
families (which causes singularities of $Quot$) but the
Euler characteristic of Ext's  is preserved under deformations
(which explains the smoothness
of $RQuot$). 

\vskip .3cm

\noindent {\bf (0.4)}
The derived $Quot$ scheme we construct is suitable for construction
of the derived moduli space of vector bundles. In the particular
case of the Hilbert scheme, i.e.,  $\Cal F=\Cal O_X$, there is
another natural derived version, $RHilb_h(X)$, which is suitable
for construction of  derived moduli spaces of algebraic varieties,
(stable) maps etc. Its construction will be carried out
in a sequel to this paper \cite{CK}. To highlight the difference
between $RHilb_h(X)$ and $RQuot_h(\Cal O_X)$, take a $\Bbb K$-point
of $Hilb_h(X)$ represented by a subscheme $Z\subset X$ with the
sheaf of ideals $\Cal J\i \Cal O_X$. Then, for a smooth $Z$ and $X$
it will be shown in \cite{CK} that
$$H^i T^\bullet_{[Z]} RHilb_h(X) =  H^{i}(Z, \Cal N_{Z/X})$$
which is  smaller than $\text{Ext}^i_X(\Cal J, \Cal O_X/\Cal J) 
= \text{Ext}^{i+1}_X(\Cal O_Z, \Cal O_Z)$ which involves the cohomology
of the higher exterior powers of the normal bundle. 

\vskip .3cm

\noindent {\bf (0.5)} The paper is organized as follows. In Section 1
we give a background treatment of the $Quot$ schemes. If we view
$Quot$ as an algebro-geometric instance of the space $Z$ from
(0.1)(A), then the role of the bigger space $C$ is played by the
ambient space of the Grassmannian embedding of $Quot$ constructed by
Grothendieck. We improve upon existing treatments 
 by
exhibiting an explicit system of  equations of $Quot$
in the product of Grassmannians (Theorem 1.4.1).

In Section 2 we make precise what we understand by the
``right derived category of schemes" in which $RQuot$ will
lie.  We develop the necessary formalism of smooth resolutions,
homotopy fiber products etc. 

In Section 3  we address a more algebraic problem: given an algebra
$A$ and a finite-dimensional $A$-module $M$,
construct the derived version of the space (called the $A$-Grassmannian)
parametrizing $A$-submodules in $M$ of the given dimension.
This construction will serve as a springboard for
constructing the derived $Quot$ scheme.

Finally, in Section 4 we give the  construction of $RQuot$,
using the approach of Section 3 and Theorem 1.4.1 which
allows us to identify $Quot$ with a version of the $A$-Grassmannian,
but for a graded module over a graded algebra. 

\vskip .3cm

\noindent {\bf (0.6)} The first published reference for the DDT program
seems to be the paper \cite{Kon} by M. Kontsevich, who  gave
an exposition of the ensuing ``hidden smoothness philosophy"
in a lecture course in Berkeley in 1994. We are also aware of earlier
unpublished suggestions of P. Deligne and V. Drinfeld containing
 very similar basic ideas. We gladly acknowledge
our intellectual debt to our predecessors. We are also grateful
to participants of the deformation theory seminar at Northwestern,
where this work originated and was reported. Both authors
were partially supported by NSF.

\vfill\eject

\heading{1. Grothendieck's Quot Scheme}\endheading

\noindent {\bf (1.1) Elementary properties.} We recall briefly the definition 
and main properties of Grothendieck's
$Quot$ scheme
(\cite{Gr}, see also \cite{Kol}\cite{Vi} for detailed treatments). Let 
$\Bbb K$ be a field and
$X$ be 
a projective scheme over $\Bbb K$, with a chosen
very ample invertible sheaf
$\Cal {O}_X(1)$. For any coherent sheaf $\Cal  G$ on $X$ denote as usual
$\Cal  G(n):=\Cal  G\otimes\Cal  O_X(n)$. The {\it Hilbert polynomial}
$h^{\Cal  G}$ of $\Cal  G$ is the polynomial in $\Bbb Q [t]$ defined by
$$h^{\Cal  G}(n)=\chi(\Cal  G(n)).$$
By Serre's vanishing theorem  $h^{\Cal  G}(n)=\operatorname{dim}
H^0(X, \Cal  G(n))$ for $n>>0$.

Now fix a coherent sheaf $\Cal  F$ on $X$, a polynomial $h'\in \Bbb
Q[t]$, and set $h:=h^{\Cal  F}-h'$. Informally the $Quot$ scheme can be 
thought of
as a ``Grassmannian of subsheaves in $\Cal  F$''; its closed
points are in 1-1 correspondence with
$$Sub_h(\Cal  F):=\{\Cal  K\subset\Cal  F\ |\ h^{\Cal  K}=h\},$$
or, equivalently, with
$$Quot_{h'}(\Cal  F):=\{\Cal  F\twoheadrightarrow\Cal  G\ |\ h^{\Cal 
G}=h'\}/Aut(\Cal  G).$$

The scheme structure reflects how quotients of $\Cal  F$ vary in 
families. More precisely, for any scheme $S$, let $\pi_X$ denote the canonical 
projection $X\times S\longrightarrow X$. Grothendieck's theorem  
then states:

\proclaim{ (1.1.1) Theorem} There exists a projective scheme $Sub_h(\Cal  F)$
(or $Quot_{h'}(\Cal  F)$) such that for any scheme $S$ we have
$$Hom(S, Sub_h(\Cal  F))=\left \{\Cal  K\subset \pi_X^*\Cal  F\ \mid \left
.\matrix \ \pi_X^*\Cal  F/\Cal  K
\ \text{is\ flat\ over\ } S\ \text{and\ has}\\ \text{ relative\ Hilbert\
polynomial}\ h'\endmatrix\ \right .\right\} .$$
\endproclaim 

Thus, in particular, we have the
 {\it universal exact sequence} on $Sub_h(\Cal  F)\times X$,
with $\Cal  S$   corresponding
to the identity map $Sub_h(\Cal  F)\to\ Sub_h(\Cal  F)$;
$$0\longrightarrow \Cal  S\longrightarrow \pi_X^*\Cal  F
\longrightarrow \Cal  Q\longrightarrow 0.\leqno (1.1.2)$$
 
\vskip 8pt

The following statement is obtained easily 
by taking $S = \op{Spec}(\Bbb K[x]/x^2)$ in Theorem 1.1.1, see
\cite{Gr, Cor. 5.3}.

\proclaim {\bf (1.1.3) Proposition} Let $[\Cal  K]$ be the $\Bbb K$-point in 
$Sub_h(\Cal  F)$ determined
by a subsheaf $\Cal  K\subset\Cal  F$ with $h^{\Cal  K}=h$. Then
the tangent space to $Sub_h(\Cal  F)$ at $[\Cal  K]$ is
$$T_{[\Cal  K]}Sub_h(\Cal  F)=\operatorname{Hom}_{\Cal  O_X}(\Cal  K,\Cal  
F/\Cal  K).$$
\endproclaim

\vskip .2cm

\noindent {\bf (1.2) The Grassmannian embedding.}
Let $W$ be a finite-dimensional vector space. By $G(k,W)$ we denote
the Grassmannian of $k$-dimensional linear subspaces
in $W$. Thus, to every such subspace $V\subset W$ there corresponds
a point $[V]\in G(k, W)$. We denote by $\widetilde V$ the tautological
vector bundle on $G(k, W)$ whose fiber over $[V]$ is $V$. 
It is well known that $T_{[V]}G(k,W)\simeq \operatorname{Hom}(V, W/V)$.

\vskip .1cm

Let $X, \Cal  O(1)$ be as before.
Set $A:=\bigoplus_{i\geq 0} H^0(X,\Cal   O_X(i))$. This is 
a finitely generated graded commutative algebra. For a coherent
sheaf $\Cal  G$ on $X$ let $\operatorname{Mod}(\Cal  G) = \bigoplus_i
H^0(X, {\Cal  G}(i))$ be the corresponding graded $A$-module. Similarly,
for a finitely generated graded $A$-module $M$ we denote by
$\text{Sh}(M)$ the coherent sheaf on $X$ corresponding to $M$
by localization. 

If $M$ is a graded $A$-module, we denote $M_{\geq p}$ the submodule
consisting of elements of degree at least $p$. Similarly, for $p\leq q$
we set $M_{[p,q]}=M_{\geq p}/M_{\geq q}$ to be the truncation of $M$
in degrees $[p,q]$. 

Given finitely generated graded $A$-modules $M, N$ we define
$$\text{Hom}_{\Cal  S}(M, N) = \lim_{\buildrel p\over\rightarrow}
\operatorname{Hom}^0(M_{\geq p},
N_{\geq p}),\leqno (1.2.1)$$
where  $\operatorname{Hom^0}$ is the
set of $A$-homomorphisms of degree $0$.

Recall the classical theorem of Serre \cite{Se, \S 59}.

\proclaim{(1.2.2) Theorem}
The category $Coh(X)$ of coherent sheaves on $X$ is equivalent to  the
category $\Cal  S$ whose objects are finitely generated graded $A$-modules 
and morphisms are  given by (1.2.1). More precisely, if $M, N$ are
objects of $\Cal  S$, then
$$\operatorname{Hom}_{\Cal  S}(M, N) = \operatorname
{Hom}_{\Cal  O_X}(\text{Sh}(M),
\text{Sh}(N)).$$
 Further, the limit in (1.2.1) is
achieved for some $p=p(M,N)$. 
\endproclaim

Part (a) of the following theorem is also due to Serre \cite{Se, \S 66}.

\proclaim{(1.2.3) Theorem} (a) For any coherent sheaf $\Cal  G$ on $X$
 there exists an
integer $p=p(\Cal  G)$ such that 
$H^j(X,\Cal  G(r))=0$ for
all $j\geq 0$ and all $r\geq p$, and
the multiplication map
$$H^0(X,\Cal  O_X(i))\otimes H^0(X,\Cal  G(r))\longrightarrow
H^0(X,\Cal  G(r+i))$$
is surjective for all $i\geq 0$ and all $r\geq p$.

(b) The number $p$  in part (a) can be chosen  uniformly
with the above properties for all subsheaves $\Cal  K$ of a fixed
coherent sheaf $\Cal  F$ on $X$ with fixed Hilbert polynomial $h^{\Cal 
K}=h$, and for all respective quotients $\Cal  F/\Cal  K$.

\endproclaim

Part (b) is proved in \cite{Mu, Lecture 14} or \cite{Vi, Thm. 1.33}.
More precisely, the discussion of \cite{Vi} is, strictly speaking,
carried out only for the
case $\Cal  F=\Cal  O_X^n$. This, however, implies the case $\Cal   F = \Cal  
O_X(i)^n$
for any $i$ and $n$ and then the case of 
an arbitrary $\Cal   F$ follows from 
this by taking a surjection $\Cal  O_X(i)^n\to\Cal  F$.

\vskip .15cm

In terms of the associated module $N=\text{Mod}(\Cal  G)$, part (a)
means that  $N_{\geq p}$ is generated by $N_p$ and
$\dim N_{r} = h^{\Cal  G}(r)$ for $r\geq p$.

\vskip .2cm

Fix now a coherent sheaf $\Cal  F$ and a polynomial $h$ and pick
$p$ such as in (1.2.3)(b) which is large enough so that the statements
of (a) hold for $\Cal  F$ as well.  
Consider
the universal exact sequence (1.1.2).  
For $r\geq p$, twisting by $\pi_X^*\Cal  O_X(r)$ and pushing forward to
$Sub_h(\Cal  F)$ produces an exact sequence of {\it vector bundles}
$$0\longrightarrow (\pi_{Sub})_*\Cal  S(r)\longrightarrow
M_r\otimes\Cal  O_{Sub}
\longrightarrow (\pi_{Sub})_*\Cal  Q(r)\longrightarrow 0,$$
with $\operatorname{rank}(\pi_{Sub})_*\Cal  S(r)=h(r)$, which in turn 
determines a map 
$$\alpha_r: Sub_h(\Cal  F)\longrightarrow G(h(r),M_r),\leqno (1.2.4)$$
Now Grothendieck's Grassmannian embedding  is as follows.

\proclaim{(1.2.5) Theorem} For $r\gg 0$ the map $\alpha_r$
 identifies  $Sub_h(\Cal  F)$ with a closed subscheme of the Grassmannian
$G(h(r),M_r)$. 
\endproclaim

\vskip .2cm

\noindent {\bf (1.3) The $A$-Grassmannian.} We now discuss a
more elementary construction which can be seen as a finite-dimensional analog
of the $Quot$ scheme.

Let  $A$ be an associative algebra over $\Bbb K$
(possibly without unit)
and $M$ be a finite-dimensional left $A$-module.  The $A$-Grassmannian is the 
closed subscheme $G_A(k,M)\subset G(k,M)$
formed by those $k$-dimensional subspaces which are left $A$-submodules.
It can be defined as the (scheme-theoretical) zero locus of the canonical
section 
$$s\in\Gamma(G(k, M), \, \Cal  Hom (A\otimes_{\Bbb K} \widetilde V, 
M/\widetilde V))
\leqno (1.3.1)$$
whose value over a point $[V]$ is the composition of the $A$-action
$A\otimes V\to M$ with the projection $M\to M/V$.
It follows that
if $V$ is  a submodule, then
$$T_{[V]}G_A(k,M) = \operatorname{ Hom}_A(V, M/V) \subset \operatorname{ 
Hom}_\Bbb K(V, M/V)=
T_{[V]}G(k,M).\leqno (1.3.2)$$
This is similar to (1.1.3).

Next, suppose that $M= \bigoplus_i M_i$ is a finite-dimensional 
$\Bbb Z$-graded
vector space,
i.e.,  each $M_i$ is finite-dimensional and $M_i=0$ for almost all $i$.
Let $k=(k_i)$ be a sequence of nonnegative integers. We denote
$G(k, M) = \prod G(k_i, M_i)$; in other words, this is the variety
of graded subspaces $V=\bigoplus V_i\subset M$ such that $\dim(V_i)=k_i$. 
As before, we denote by $[V]$ the point of $G(k, M)$ represented
by a graded subspace $V$, and denote by $\widetilde V = \bigoplus\widetilde 
V_i$
the tautological graded vector bundle over $G(k, M)$.

Let now $A = \bigoplus_i A_i $ be a $\Bbb Z$-graded associative algebra
and $M = \bigoplus_i M_i$ be a finite-dimensional graded left $A$-module,
  We have then the graded $A$-Grassmannian $G_A(k, M)\subset G(k, M)$ 
parametrizing graded $A$-submodules $V\subset M$.
It can be defined as the common zero locus of the natural sections
$s_{ij}$ of the bundles $\Cal  Hom (A_i\otimes\widetilde V_j,
M_{i+j}/\widetilde V_{i+j})$. 
For  a submodule
$V$ we have
$$T_{[V]} G_A(k, M) = \operatorname{Hom}_A^0 (V, M/V),\leqno (1.3.3)$$
where $\operatorname{Hom}_A^0$ means the set of homomorphisms of degree 0. 

\vskip .3cm

\noindent {\bf (1.4)  $Quot$ as an $A$-Grassmannian.} We specialize
the considerations of (1.3) to the case 
 $$A = \bigoplus H^0(X, \Cal  O_X(i)), \quad M=\bigoplus_i H^0(X, \Cal  F(i)),
\quad \Cal  F\in \text{Coh}(X) $$ 
from (1.2). Let $p>0$ be chosen as in (1.2). For $p\leq q$ the morphism
$$\alpha_{[p,q]} =\prod_{r=p}^q \alpha_r: Sub_h(\Cal  F)\to
\prod_{r=p}^q G(h(r), M_r) = G(h, M_{[p,q]})$$
takes values, by construction, in the $A$-Grassmannian $G_A(h, M_{[p,q]})$.
The following result extends Theorem 1.2.5 by providing explicit
relations for the Grassmannian embedding of $Quot$.
It seems not to be found in the literature.
A related statement (which does not take into account the nilpotents
in the structure sheaves of the schemes involved and assumes
$\Cal F =\Cal O_X$), is due to Gotzmann [Go, Bemerkung 3.3]. In his situation
it is enough to take $q=p+1$. 

\proclaim{(1.4.1) Theorem} For $0\ll p\ll q$ the morphism
$\alpha_{[p,q]}: Sub_h(\Cal  F)\to G_A(h, M_{[p,q]})$ is an isomorphism. 
\endproclaim

First notice that we may assume that $\Bbb K$ is algebraically closed.
 Before giving the proof of the theorem, we need some preparations. 
To unburden the notation, for $q\geq p$ set
$$G_q:=G_A(h,M_{[p,q]}).$$
In particular $G_p=G(h(p), M_p)$, since $M_{[p,p]}=M_p$ has 
trivial $A$-module structure.

For $r\geq s\geq p$, let $\varphi_{rs}:G_r\longrightarrow G_s$ be the 
canonical projection. We have then an {\it inverse system} of schemes
$$\minCDarrowwidth{7 mm}
\CD G_p@<<< G_{p+1}@<<<\cdots @<<< G_r @<<< \cdots .\endCD\leqno (1.4.2)$$

\proclaim{(1.4.3) Lemma} The projective limit of the diagram (1.4.2) 
in the category of schemes exists and is identified with $Sub_h(\Cal  F)$.
\endproclaim

\demo{Proof} We have to show that for any scheme $S$, a compatible 
system of maps $S\longrightarrow G_r$, $r\geq p$ gives rise to a map
$S\longrightarrow Sub_h(\Cal  F)$. But such a system 
gives rise to a family (parametrized by $S$) of graded $A$-submodules
of $M_{\geq p}$ with Hilbert polynomial 
(in fact, even the Hilbert function) equal to $h$, i.e., to a graded
$A\otimes \Cal  O_S$-submodule $\Cal  V\subset M\otimes \Cal  O_S$ such that
each $\Cal  V_i$ is a  projective $\Cal  O_S$-module (i.e.,
a vector bundle on $S$) of rank $h(i)$. It follows
that each $(M_i\otimes\Cal  O_S)/\Cal  V_i$ 
is projective, as an $\Cal  O_S$-module. Thus the graded $\Cal  O_S$-module
$(M\otimes \Cal  O_S)/\Cal  V$ is flat (because it is
the union  of projective
$\Cal  O_S$-submodules $(M\otimes \Cal  O_S)_{\leq p}/\Cal  V_{\leq p}$).
By Serre's theorem $\Cal  V$ gives a subsheaf $\Cal  K\subset \pi^*\Cal  F$
and the quotient sheaf  $\pi^*\Cal  F/\Cal  K$, corresponding to
the $\Cal  O_S$-flat graded $A\otimes\Cal  O_S$-module 
$(M\otimes \Cal  O_S)/\Cal  V$, is also flat over $\Cal  O_S$.
Hence we get the required map
$S\longrightarrow Sub_h(\Cal  F)$.
\enddemo

Recall that  any morphism $f: Y\to Z$ of projective schemes
has a well-defined image which is a closed subscheme $\text{Im}(f)\subset
Z$ satisfying the usual categorical universal property. 

With this understanding, for any $r\geq p$,  we consider the subscheme 
$\widetilde{G}_r$ of $G_r$
defined by 
$$\widetilde{G}_r:=\bigcap_{r'\geq r}\operatorname{Im}
\{\varphi_{r'r}:G_{r'}\longrightarrow G_r\}.\leqno (1.4.4)$$
Because of the Noetherian property, the intersections in (1.4.4)
in fact stabilize.

\proclaim{(1.4.5) Lemma} Together with the restrictions of the natural
projections, the subschemes $\widetilde{G}_r$ form an inverse system
of surjective maps with the same projective limit $Sub_h(\Cal  F)$ as
the system (1.4.2).
\endproclaim 

\demo{Proof}
This is a purely formal argument. We consider $(G_r)$ and
$(\widetilde G_r)$ as pro-objects in the category of schemes
(see \cite{GV}, \S 8)
and will show that they are isomorphic in the category of
pro-objects. This will imply that $\lim\limits_\leftarrow \widetilde G_r$
exists and is isomorphic to $\lim\limits_\leftarrow  G_r$. 
First, the componentwise morphism
of inverse systems $(u_r: \widetilde G_r \to G_r)$ gives
a morphism of pro-objects, which we denote $u_*$. Next,
  stabilization of the images implies
that for every $r$ there is a $q=q(r)$ and a morphism $v_r: G_{q(r)}\to 
\widetilde G_r$. These constitute a morphism of pro-objects
$v_*: (G_r)\to (\widetilde G_r)$,  which one checks is inverse to $u_*$. 
\enddemo

\proclaim{(1.4.6) Corollary} The projective
system $(\widetilde G_r)$ is constant. In particular, for any $r\geq p$
  the natural
projection $\widetilde \varphi_r: Sub_h(\Cal  F)\to \widetilde G_r$
is an isomorphism.
\endproclaim

\demo{Proof} This follows from the previous two lemmas and Grothendieck's
theorem 1.2.5 on the Grassmannian embedding which can be formulated
by saying
that $\widetilde\varphi_p: Sub_h(\Cal  F)\to \widetilde G_p$
is an isomorphism. 
 \enddemo

 \demo{Proof of Theorem (1.4.1)}
It follows from Lemmas 1.4.3 and 1.4.5  and Corollary 1.4.6
that we have a commutative diagram
$$\matrix  & & & & & & Sub_h(\Cal  F) \\
  & & & & & & \downarrow \\
\widetilde{G}_p &
\buildrel\sim\over\longleftarrow & \widetilde{G}_{p+1} &
\buildrel\sim\over\longleftarrow & \cdots & \buildrel\sim\over\longleftarrow & 
\displaystyle{\lim_{\longleftarrow}}\;\widetilde{G}_{r} \\
 \bigcap & & \bigcap & & & & || \\
G_p & \longleftarrow & G_{p+1} & \longleftarrow & \cdots &
\longleftarrow &\displaystyle{\lim_{\longleftarrow}}\; G_{r} \\
\endmatrix $$
such that for every $r\geq p$ the induced map
$Sub_h(\Cal  F)\rightarrow\widetilde{G}_r$ is an isomorphism
and $\alpha_{[p,r]}$ factors as 
$Sub_h(\Cal  F)\rightarrow\widetilde{G}_r\hookrightarrow G_r$.
Let $q$ be such that for any $r\geq q$
$$\operatorname{Im}\{\varphi_{rp}:G_{r}
\longrightarrow G_p\}=\widetilde{G}_p(\cong Sub_h(\Cal  F)).$$ 
 
\vskip 10pt

For an $A$-submodule
$V=\bigoplus_{i=p}^rV_i$ of
$M_{[p,r]}$ 
let $\vec{V}$ be the $A$-submodule of $M_{\geq p}$ generated
by $V$. Let now $W=W_p\in \widetilde{G}_p$ be a subspace of $M_p$.
Since $\widetilde{G}_p\cong Sub_h(\Cal  F)$ we have $W=H^0(X,\Cal  K(p))$,
for some $\Cal  K\subset\Cal  F$ with $h^{\Cal  K}=h$. By our choice of $p$, 
it follows that $\vec{W}=\bigoplus_{i \geq p}H^0(X,\Cal  K(i))$. 

\proclaim{(1.4.7) Lemma} For $r\geq q$ the set-theoretic fibre 
$\varphi_{rp}^{-1}(W)\subset G_r$ consists of the unique $\Bbb K$-point  
$(\vec{W})_{[p,r]}:=\vec{W}\cap M_{[p,r]}$.
\endproclaim

Granting this for a moment, it follows that for $r\geq q$ the map
$\varphi_{rp}$ is a bijection on $\Bbb K$-points onto
$\widetilde{G}_p$, and therefore
$\alpha_{[p,r]}$ is also a bijection on $\Bbb K$-points onto
$G_r=G_A(h,M_{[p,r]})$. This gives  Theorem 1.4.1 at the level of sets.
To prove it in general, consider the tautological family of $A$-submodules
of $M_{[p,r]}$ over $G_A(h,M_{[p,r]})$ obtained by restricting the 
tautological vector subbundle over $G(h,M_{[p,r]})$. This determines (by
pull-back to $G_A(h,M_{[p,r]})\times X$ and application of the functor
``\; $\vec{*}$\; '') a family of $A$-submodules of $M_{\geq p}$ with Hilbert 
polynomial equal to $h$. The same argument as in the proof of Lemma 1.4.3
gives then a map $\beta_{[p,r]}:G_A(h,M_{[p,r]})\longrightarrow
Sub_h(\Cal  F)$ which is easily seen to be an inverse for $\alpha_{[p,r]}$.
\enddemo

\demo{Proof of Lemma 1.4.7}
Let $V\in G_A(h,M_{[p,r]})$ be such that $V\cap M_p=W$. Then 
$\vec{W}\subset \vec{V}\subset M_{\geq p}$. Since $V$ is an $A$-submodule, 
it follows that
$(\vec{W})_{[p,r]}\subset V$. But for each $p\leq i\leq r$ the dimension 
over $\Bbb K$ of the graded components of degree $i$ of this last two modules
is the same, therefore $(\vec{W})_{[p,r]}= V$.
\enddemo

This concludes the proof of Theorem 1.4.1.

\vfill\eject

\heading {2. The right derived category of schemes.}\endheading

In this section we develop the minimal necessary background
suitable for taking right derived functors on the category
of schemes (which correspond to left derived functors
on the category of commutative algebras). 

\vskip .3cm

\noindent {\bf (2.1) Dg-vector spaces, algebras
and modules.} From now on we assume\footnote
{The reasons for the characteristic 0 assumptions
are the standard ones in the theory of dg-algebras
\cite{Le} \cite{Q2}. For example,
the construction of M-homotopies in (3.6)
requires taking anti-derivatives of polynomials with
coefficients in $\Bbb K$.}
that the base field  $\Bbb K$ 
has characteristic 0. By a complex (or dg-vector space)
we always mean a cochain complex, i.e., a graded vector space
$C$ with a differential of degree $+1$. The grading
of  complexes will be always indicated in the superscript,
to distinguish it from other types of grading which may be eventually
present (such as in (1.2) above). If $C$ is a complex and $a\in C^i$,
 we write $\bar a=i$. We also write $H(C)$ for the graded space
of cohomology of $C$ and $C_\# $ for the graded vector
space obtained from $C$ by forgetting the differential. 
A morphism $f: C\to D$ of complexes is called a quasiisomorphism
if $H(f): H(C)\to H(D)$ is an isomorphism. 

Complexes form a symmetric monoidal category $dgVect$ with respect to the
usual tensor product and the symmetry operator given by the Koszul
sign rule: $a\otimes b\mapsto (-1)^{\bar a\bar b} b\otimes a$. 
By an associative resp. commutative  dg-algebra we mean
an associative, resp. commutative algebra in $dgVect$. By a graded
algebra we mean a dg-algebra with zero differential. 
Thus for every dg-algebra $A$ we have graded algebras $H(A), A_\# $.
Note that a graded commutative algebra in this sense satisfies
$ab=(-1)^{\bar a\bar b} ba$. 


Similar conventions and terminology will be used for dg-modules 
over a dg-algebra
$A$ (left or right, if $A$ is not commutative).

In this paper we will always consider (unless otherwise specified),
only dg-algebras $A$ which are $\Bbb Z_-$-graded, i.e., have $A^i=0$ for
$i>0$.

The following remark, though obvious, is crucial for gluing commutative
dg-algebras into more global objects.

\proclaim{(2.1.1) Proposition} If $A$ is a $\Bbb Z_-$-graded associative
dg-algebra and $M$ is a left dg-module over $A$, then each
$d: M^i\to M^{i+1}$ is $A^0$-linear.
\endproclaim

Let $A$ be an associative dg-algebra. A left dg-module $M$ over $A$
is called quasifree, if $M_\#$ is free over $A_\#$, so as a graded
module, $M= A\otimes_{\Bbb K} E^\bullet$, where $E^\bullet$
is some graded vector space of generators.

Let  $M, N$ be left dg-modules
over $A$. Morphisms of dg-modules $M\to N$ are degree 0 cocycles
in the cochain complex $\op{Hom}^\bullet_{A}(M, N)$ which
consists of all $A$-linear morphisms and whose differential
is given by the commutation with the differentials in $M, N$. 
Two morphisms $f, g: M\to N$ are called homotopic, if they are cohomologous
as cocycles, i.e., if there exists a morphism $s: M\to N[-1]$ of
$A_\#$-modules such that $d_N s+sd_M = f-g$. 

\proclaim{(2.1.2) Proposition} Let $M, N$ be two $\Bbb Z_-$-graded
left dg-modules over $A$, and suppose that $M$ is quasifree
and $N$ is acyclic in degrees $<0$. Let $f, g: M\to N$ be dg-morphisms
which induce the same morphism $H^0(M)\to H^0(N)$. Then $f$ is homotopic
to $g$.
\endproclaim

\noindent {\sl Proof:} This is  a standard inductive
construction and is left to the reader. In fact, in Proposition
3.6.4 we give a less trivial, nonlinear version of this construction
and spell out the proof in full detail. The reader can
easily adapt that proof to the present linear situation. 

\vskip .1cm

\proclaim{(2.1.3) Corollary} If $M$ is a quasifree dg-module over
$A$ which is bounded from above and which is exact with respect
to its differential, then $M$ is contractible, i.e., its
identity map is homotopic to 0.
\endproclaim

\noindent {\sl Proof:} By shifting the degree we can assume that
$M$ is $\Bbb Z_-$-graded. Then apply (2.1.2) to $M=N$ and to $f= \op{Id}_M$,
$g=0$.

\vskip .1cm

\proclaim {(2.1.4) Corollary} A quasiisomorphism $f: M\to N$ of
quasifree, bounded from above dg-modules over $A$, is
a homotopy equivalence.
\endproclaim

\noindent {\sl Proof:} Consider the dg-module $\op{Cone}(f)$. It is
acyclic, quasifree and bounded from above, so contractible
by (2.1.3). This implies our statement.

\vskip .1cm

\proclaim {(2.1.5) Corollary} In the situation of Corollary 2.1.4,
the dual morphism
$$f^*: \op{Hom}_A(N, A)\to \op{Hom}_A(M, A)$$
is a quasiisomorphism.
\endproclaim

\noindent {\sl Proof:} 
As $f$ is a homotopy equivalence, so is $f^*$, because 
homotopies are inherited under functorial constructions
on modules such as $\op{Hom}_A(-, A)$.

\vskip .3cm

\noindent {\bf (2.2) Dg-schemes.}

\proclaim{(2.2.1) Definition} (a) By  a dg-scheme we mean a pair
$X=(X^0, \Cal  O_X^\bullet)$, where $X^0$ is an ordinary scheme
and $\Cal  O_X^\bullet$ is a sheaf of ($\Bbb Z_-$-graded)
commutative dg-algebras on $X^0$ such that $\Cal  O_X^0 = \Cal  O_{X^0}$
and each $\Cal  O_X^i$ is quasicoherent over $\Cal  O_X^0$. 

(b) A morphism $f: X\to Y$ of dg-schemes is just a morphism of (dg-) ringed
spaces, i.e., a morphism $f_0: X^0\to Y^0$ of schemes together with a morphism
of sheaves of dg-algebras $f_0^*{\Cal  O}_Y^\bullet\to {\Cal  O}_X^\bullet$.
The category of dg-schemes will be denoted by $dgSch$.

\endproclaim

By a graded scheme we mean a dg-scheme $X$ in which $\Cal  O_X^\bullet$
has trivial differential. Any ordinary scheme will
be considered as a dg-scheme with trivial grading and differential. 

By (2.1.1), for a dg-scheme $X$ each $d: \Cal  O_X^i\to\Cal  O_X^{i+1}$
is $\Cal  O_{X^0}$-linear and hence $\underline{H}^i (\Cal  O_X^\bullet)$
are quasicoherent sheaves on $X^0$. We define the ``degree 0 truncation"
of $X$ to be the ordinary scheme
$$\pi_0(X) = \text{Spec} \underline{H}^0(\Cal  O_X^\bullet)
\subset X^0. \leqno (2.2.2)$$
The notation is chosen to suggest analogy with homotopy groups
in topology. Note that for any ordinary scheme $S$ we have
$$\text{Hom}_{dgSch}(S, X) = \text{Hom}_{Sch}(S, \pi_0(X)).\leqno (2.2.3)$$
 Note also that each $\underline{H}^i(\Cal  O_X^\bullet)$
can be regarded as a quasicoherent sheaf on $\pi_0(X)$. 
We have then two graded schemes naturally associated to $X$:
$$X_\#  = (X^0, \Cal  O^\bullet_{X, \# }), \quad X_h =
 (\pi_0(X), \underline{H}^\bullet (\Cal  O_X^\bullet)).\leqno (2.2.4)$$
A morphism $f: X\to Y$ of dg-schemes will be called a quasiisomorphism
if the induced morphism of graded schemes $f_h: X_h\to Y_h$
is an isomorphism. 
We denote by $\Cal  DSch$ the category obtained from $dgSch$
by inverting all quasiisomorphisms and call it the (right)
derived category of
schemes. It is suitable for taking right derived functors on
 schemes (which correspond to left derived functors at the level 
of commutative algebras).

\vskip .2cm

\noindent {\bf (2.2.5) Examples.} (a)  If $A^\bullet$ is a commutative
$\Bbb Z_-$-graded 
dg-algebra, we have a dg-scheme $X = \text{Spec}(A^\bullet)$ defined as 
follows.
The scheme $X^0$ is $\text{Spec}(A^0)$, and the sheaf ${\Cal  O}_X^i$ is
the quasicoherent sheaf on $\text{Spec}(A^0)$ associated 
to the $A^0$-module $A^i$.
A dg-scheme $X$ having the form $\text{Spec}(A^\bullet)$ will be called affine.
We will also write $A^\bullet = {\Bbb K}[X]$ and call $A$ the coordinate (dg-) 
algebra of $X$.

(b) In particular, if $E^\bullet$ is a $\Bbb Z_+$-graded complex
of finite-dimensional vector spaces, we have the dg-scheme
$|E|=\operatorname{ Spec}\, S(E^*)$, which is ``the linear dg-space
$E^\bullet$ considered as a scheme". We will extend this notation
as follows. Let   $F^\bullet$ be a $\Bbb Z_-$-graded complex of
possibly infinite-dimensional vector spaces. Then we write
$|F^*| =  \operatorname{ Spec}\, S(F)$. 
We will use this notation for ungraded vector spaces as well. 

\vskip .2cm

  \proclaim{ (2.2.6) Definition} A dg-scheme $X$ is said to be of finite type,
if $X^0$ is a scheme of finite type, and  each $\Cal  O_X^i$ is coherent,
as a sheaf on $X^0$. We say that $X$ is (quasi)projective if it is
of finite type and $X^0$ is (quasi)projective in the usual sense. 
\endproclaim

\vskip .2cm

\noindent {\bf (2.3) Dg-sheaves.} We now globalize the usual theory of 
dg-modules
over a dg-algebra \cite {HMS} \cite{Mc}.

\proclaim {(2.3.1) Definition} A quasicoherent (dg)-sheaf on a dg-scheme $X$ 
is a sheaf
${\Cal  F}^\bullet$ of ${\Cal  O}_{X}^\bullet$-dg-modules such that every
${\Cal  F}^i$ is quasicoherent over ${\Cal  O}_{X^0}$. If $X$ is of finite 
type,
we shall say that a quasicoherent sheaf ${\Cal  F}^\bullet$ is coherent, if 
each ${\Cal  F}^i$
is coherent over ${\Cal  O}_{X^0}$
and if $\Cal F^\bullet$ is bounded above, i.e., $\Cal  F^i=0$ for $i\gg 0$.  
\endproclaim

If $\Cal  F$ is a quasicoherent dg-sheaf on $X$, we have
graded sheaves
 $\Cal  F_\#$ on $X_\#$ and $\underline{H}^\bullet(\Cal  F^\bullet)$
on $X_h$.

Morphisms  of quasicoherent dg-sheaves $\Cal F^\bullet\to\Cal G^\bullet$
and homotopies between such morphisms
are defined in the obvious way, cf. (2.1). 
 
A morphism $\Cal F^\bullet\to\Cal G^\bullet$ is called a quasiisomorphism,
if  the induced morphism of graded sheaves $\underline{H}^\bullet(\Cal  
F^\bullet)\to \underline{H}^\bullet(\Cal  G^\bullet)$ is an isomorphism.

We denote by $\Cal{DQC}oh_X$ (resp. $ \Cal {DQC}oh^-_X$) the derived category
of quasicoherent dg-sheaves (resp. of bounded above quasicoherent dg-sheaves)
on $X$. Its objects are dg-sheaves of the described kind
and morphisms are obtained
 by first passing to
homotopy classes of morphisms and then localizing the resulting category
by quasiisomorphisms. Similarly, if $X$ is of finite type, we
have $\Cal {DC}oh_X$, the derived category of coherent dg-sheaves.
These are triangulated categories naturally associated to $X$.

\vskip .1cm

If $S$ is an ordinary scheme, by a graded vector bundle we mean
a graded sheaf $E^\bullet$ of $\Cal  O_S$-modules such that
each $E^i$ is locally free of finite rank.

\proclaim {(2.3.2) Definition} Let $X$ be a dg-scheme of finite type.
A dg-vector bundle on $X$
is a coherent dg-sheaf $\Cal  F^\bullet$ such that locally, on the Zariski
topology of $X^0$, the graded sheaf 
$\Cal  F^\bullet_\#$ of $\Cal  O_{X,\#}^\bullet$-modules is isomorphic to
$\Cal  O_{X,\#}^\bullet\otimes_{\Cal  O_{X^0}} E^\bullet$
where $E^\bullet$ is a graded vector bundle on $X^0$ bounded from above.
\endproclaim

If $X^0$ is connected, the sequence $r=\{r_i\} = \{\text{rk}(E^i)\}$ is 
uniquely defined by $\Cal  F^
\bullet$ and is called the graded rank of $\Cal  F^\bullet$. We say that
$\Cal  F^\bullet$ has bounded rank, if $r_i=0$ for $i\ll 0$.
In this case $\Cal  F^* = \Cal  Hom_{\Cal  O_X^\bullet}
(\Cal  F^\bullet, \Cal  O_X^\bullet)$
is again a dg-vector bundle of bounded rank.

\noindent {\bf (2.3.3) Notation.} If $\Cal  A$ is a quasicoherent
sheaf of $\Bbb Z_-$-graded dg-algebras on a dg-scheme $X$, then
we have a dg-scheme $\text{Spec}(\Cal  A)\to X$. 
If $\Cal  F^\bullet$ is a dg-vector bundle on $X$ with $r_i(\Cal  F^\bullet)=0$
for $i<0$, then the symmetric algebra $S(\Cal  F^*)$ is $\Bbb Z_-$-graded
and we denote $|\Cal  F^\bullet| = \text{Spec}(S(\Cal  F^*))$. 
Similarly, if $\Cal  F^\bullet$ is a dg-vector bundle with $r_i=0$ for $i>0$, 
we write $|\Cal  F^*| = \text{Spec}(S(\Cal  F^\bullet))$.

\vskip .2cm

We now establish the existence of good resolutions of dg-sheaves
by vector bundles. The following fact is well known
(part (a) is in fact true for any scheme $S$ which can be
embedded as a closed subscheme into a smooth algebraic scheme,
see \cite{Fu} \S B.8, and part (b) follows from it). 

\proclaim {(2.3.4) Lemma} Let $S$ be a quasiprojective scheme.
Then:

(a) For any coherent sheaf $\Cal F$ on $S$ there exists a vector bundle
$E$ and a surjection $E\to\Cal F$.

(b) For any quasicoherent sheaf $\Cal F$ on $S$ there exists an $\Cal O_S$-flat
quasicoherent sheaf $\Cal E$ and a surjection $\Cal E\to\Cal F$. 

\endproclaim

\proclaim{(2.3.5) Proposition} Let $X$ be a quasiprojective dg-scheme
and let $\Cal F^\bullet$ be a quasicoherent dg-sheaf
on $X$ which is bounded from above.
Then:

(a) $\Cal F^\bullet$ is quasiisomorphic to a quasicoherent dg-sheaf
$\Cal E^\bullet$, still bounded above and such that $\Cal E_\#$ is
flat over $\Cal O_{X \#}$.
 
(b)  If $\Cal F^\bullet$ is coherent, it is quasiisomorphic to a vector
bundle.

\endproclaim

\demo{Proof} Standard inductive argument, using, at each step,
Lemma 2.3.4 and a dg-module of the form $\Cal O_X^\bullet\otimes_{\Cal O_{X^0}}
\Cal E$ where $\Cal E$ is a flat $\Cal O_{X^0}$-sheaf or a vector
bundle, cf. \cite{Mc}, \S 7.1.1.
\enddemo

\vskip .3cm

\noindent {\bf (2.4) Derived tensor product.}
Let $X$ be a dg-scheme. 
As $\Cal O_X^\bullet$ is a sheaf of commutative dg-algebras,
we can form the tensor product $\Cal F^\bullet\otimes_{\Cal O^\bullet_X} \Cal 
G^\bullet$ of any quasicoherent dg-sheaves. We will need the derived
functor of the tensor product as well.

\proclaim {(2.4.1) Proposition} Let $\Cal E^\bullet, \Cal G^\bullet$
be  quasicoherent, bounded from above, dg-sheaves on a dg-scheme $X$ 
and suppose that $\Cal E_\#$ is flat over $\Cal O_{X\#}$. Then we have the 
converging (Eilenberg-Moore)
spectral sequence of sheaves on $X^0$
$$E_2 =  \Cal Tor_\bullet^{\underline H^\bullet(\Cal O_X)}
(\underline{H}^\bullet(\Cal  E^\bullet), \,  \underline{H}^\bullet(\Cal  
G^\bullet)) \quad \Rightarrow \quad \underline {H}^\bullet
(\Cal E^\bullet\otimes_{\Cal O^\bullet_X} \Cal G^\bullet).$$
\endproclaim

Assuming this proposition, we can make the following definition.

\proclaim {(2.4.2) Definition} Let $X$ be a quasiprojective
dg-scheme and $\Cal F^\bullet, \Cal G^\bullet$ be quasicoherent
dg-sheaves on $X$ bounded above. The derived tensor product
$\Cal F^\bullet\otimes^L_{\Cal O_X^\bullet}\Cal G^\bullet$ is, by
definition, the usual tensor product $\Cal E^\bullet\otimes 
_{\Cal O_X^\bullet}\Cal G^\bullet$ where $\Cal E^\bullet$ is any
resolution of $\Cal F^\bullet$ which is bounded above and $\#$-flat.
\endproclaim

The existence of the resolution is given by (2.3.5), the independence
of the resolution by (2.4.1).

\proclaim {(2.4.3) Proposition} In the situation of (2.4.2), we have
the two converging Eilenberg-Moore spectral sequences
$$E_1 = \Cal Tor_\bullet^{\Cal O_{X\#}}(\Cal F_\#, \Cal G_\#)\quad
\Rightarrow 
\quad \underline{H}^\bullet(\Cal F^\bullet\otimes^L_{\Cal O^\bullet_X}
 \Cal G^\bullet),$$
$$E_2 = \Cal Tor_\bullet^{\underline H^\bullet(\Cal O_X)}
(\underline{H}^\bullet(\Cal  F^\bullet), \,  \underline{H}^\bullet(\Cal  
G^\bullet) \quad \Rightarrow \quad \underline {H}^\bullet
(\Cal F^\bullet\otimes^L_{\Cal O^\bullet_X} \Cal G^\bullet).$$\endproclaim

\vskip .1cm

Propositions 2.4.1-3 are obtained by globalizing the known statements
about dg-modules over a dg-algebra, see, e.g., \cite{Mc} \S 7.1.1.
Compared to {\it loc. cit.} however, our class of algebras is more
restricted and our class of resolutions is more general, so we indicate
the main steps.

 Let $A$ be a $\Bbb Z_-$-graded commutative dg-algebra
and $P,Q$ be two dg-modules over $A$ bounded above. Then we have an
{\it ad hoc} definition of the derived tensor product based on the
bar-resolution
$$\text{Bar}_A(P) = \biggl\{ ...\to A\otimes_{\Bbb K} A\otimes_{\Bbb K} P
\to  A\otimes_{\Bbb K}P \biggr\} \buildrel qis\over\longrightarrow P.
\leqno (2.4.4)$$
 More precisely, $\text{Bar}_A(P)$ is the total complex of the
double complex inside the braces; denote by $\text{Bar}^i_A(P)$, $i\leq 0$
the $i$th column of this double complex, i.e., $A^{\otimes (i+1)}\otimes P$. 
This resolution satisfies the following properties:

\proclaim {(2.4.5) Proposition} (a) $\operatorname{Bar}_A(P)_\#$ is a free
$A_\#$-module.

(b) If $P\to P'$ is a quasiisomorphism of dg-modules bounded from above,
then $\operatorname{Bar}_A(P)\to\operatorname{Bar}_A(P')$ is a 
quasiisomorphism.

(c) $H^\bullet(\operatorname{Bar}^i_A(M)) = \operatorname 
{Bar}^i_{H^\bullet(A)}
(H^\bullet(M))$.

\endproclaim

\demo{Proof} (a) is obvious; (b) follows by a spectral sequence argument
(legitimate since the complexes are bounded from above) and
(c) follows from the K\"unneth formula.
\enddemo

We define the ad hoc derived tensor product to be
$$P\boxtimes_A Q = \text{Bar}_A(P)\otimes_A Q = 
\biggl\{ ...\to P\otimes_{\Bbb K} A \otimes_{\Bbb K} Q
\to P\otimes _{\Bbb K}Q\biggr\}.\leqno (2.4.6)$$
As before, this is in fact
 really the total complex of a double complex whose vertical
differential is induced by $d_P, d_A, d_Q$ and the horizontal one
by the structure of $A$-modules on $P,Q$. So the standard
spectral sequence of this double complex gives us the first Eilenberg-Moore
spectral sequence in the form
$$E_1 = \operatorname{Tor}_\bullet^{A_\#}(P_\#, Q_\#) \quad\Rightarrow
\quad H^\bullet (P\boxtimes_A Q).\leqno (2.4.7)$$

 \proclaim {(2.4.8) Corollary} If $F\to P$ is a quasiisomorphism with
$F$ bounded above and $F_\#$ flat over $A_\#$, then 
$P\boxtimes_A Q$ is quasiisomorphic to the usual tensor product 
$F\otimes_A Q$. In particular, $F\otimes_A Q$ is independent on
the choice of a $\#$-flat resolution $F\to P$.  
\endproclaim 

Part (c) of Proposition 2.4.5 gives the second Eilenberg-Moore spectral
sequence in the form
$$E_2 = \operatorname{Tor}_\bullet^{H^\bullet(A)}(H^\bullet(P), H^\bullet(Q))
\quad\Rightarrow \quad H^\bullet(P\boxtimes_A Q). \leqno (2.4.9)$$

\proclaim {(2.4.10) Corollary} If $P_\#$ is flat over $A_\#$, then we have
a spectral sequence converging to the ordinary tensor product
$$E_2 = \operatorname{Tor}_\bullet^{H^\bullet(A)}(H^\bullet(P), H^\bullet(Q))
\quad\Rightarrow \quad H^\bullet(P\otimes_A Q).$$
\endproclaim

Proposition 2.4.1 follows from (2.4.10) by gluing together the
spectral sequences corresponding to $A=\Gamma(U, \Cal O_X^\bullet)$
for affine open $U\i X^0$. At the same time, we get the second spectral
sequence in (2.4.3). As for the first spectral sequence, it
follows by gluing together the spectral sequence obtained 
from (2.4.7), (2.4.8) and the definition of $\otimes^L$.

\vskip .3cm

\noindent {\bf (2.5) Dg-manifolds and tangent complexes.} 

\proclaim{(2.5.1) Definition} A dg-scheme $M$ is called smooth
(or a dg-manifold) if the following conditions hold:

(a) $M^0$ is a smooth algebraic variety.

(b) Locally on Zariski topology of $M^0$, we have an isomorphism
of graded algebras
$$\Cal  O_{M^\#}^\bullet \quad\simeq \quad S_{\Cal  O_{M^0}}\biggl(Q^{-1}
\oplus Q^{-2} \oplus ...\biggr),$$
where $Q^{-i}$ are vector bundles (of finite rank) on $M^0$. 

\endproclaim

An equivalent way of expressing (b) is that every truncation
$\Cal  O_{M^\#}^{\geq -d}$ is locally, on   $M^0$
isomorphic to the similar truncation of a free graded-commutative
$\Cal  O_{M^0}$-algebra with finitely many generators in each degree. 

\vskip .1cm

The graded vector bundle $Q^\bullet_M = \bigoplus_{i\leq -1} Q^i_M$
from (b) can be defined globally as the bundle of primitive elements:
$$Q^\bullet_M = \Cal  O_M^{\leq -1} \bigl/ (\Cal  O_M^{\leq -1})^2,
\leqno (2.5.2)$$
but we do not have a natural embedding $Q^\bullet_M\to \Cal  O^\bullet_M$.

\vskip .1cm

The {\it dimension} of a dg-manifold $M$ is the sequence $\dim(M) = 
\{d_i(M)\}_{i\geq 0}$ where
$$d_0(M) = \dim(M^0), \quad d_i(M)=\text{rk}(Q^{-i}_M), i>0.\leqno (2.5.3)$$

The {\it cotangent dg-sheaf} $\Omega^{1\bullet}_M$ of $M$ is defined, as
in the commutative case, as the target of the universal derivation
 $\delta: \Cal O^\bullet_M
\to \Omega^{1\bullet}_M$. The proof of the following proposition
is standard. 

\proclaim {(2.5.4) Proposition} $\Omega^{1\bullet}_M$ is $\Bbb Z_-$-graded
coherent sheaf, which is a vector bundle of rank
 $\{r_i = d_{-i}(M)\}_{i\leq 0}$.
\endproclaim

We have the {\it tangent dg-sheaf} $T^\bullet M$ defined as usual
via derivations
$$T^\bullet M = \Cal Der_{\Bbb K} (\Cal O_M^\bullet, \Cal O_M^\bullet) = 
\Cal Hom_{\Cal O_M^\bullet}(\Omega^{1\bullet}_M, \Cal O_M^\bullet).\leqno
(2.5.5)$$
This is a quasicoherent sheaf of dg-Lie algebras on $M$.
Its differential is given by the commutator with the differential in
$\Cal O_M^\bullet$.
It is coherent if and only if $d_i(M)=0$ for $i\gg 0$. 
Further if 
$x\in M(\Bbb K) = \pi_0(M)(\Bbb K)$, is a $\Bbb K$-point, the {\it tangent
dg-space} to $M$ at $x$ is defined by
$$T^\bullet_xM = \operatorname{Der}_{\Bbb K}(\Cal O_M^\bullet, \Bbb K_x) = 
T^\bullet M \otimes_{\Cal O_M^\bullet} \Bbb K_x. \leqno (2.5.6)$$
Here  $\Bbb K_x$ is the 1-dimensional $\Cal  O_M^\bullet$-dg-module
corresponding to $x$. Similarly, we can define the tangent space
at any $\Bbb F$-point of $\pi_0(M)$, where $\Bbb K\subset\Bbb F$ is a field extension.
We will sometimes use the following suggestive ``topological" notation:
$$\pi_i(M,x) = H^{-i}(T^\bullet_x M), i<0.\leqno (2.5.7)$$
One justification of it is given by the following remark.

\proclaim{(2.5.8) Proposition} Let $\Bbb F$ be an extension
of $\Bbb K$ and $x$ be an $\Bbb F$-point of $x$. Any choice of a formal coordinate
system on $M$ near $x$ gives rise to a structure
of a homotopy Lie algebra \cite{St2} on the shifted tangent
dg-space $T^\bullet_x M[-1]$. In particular, at the level
of cohomology we have well defined ``Whitehead products"
$$[- , -]: \pi_i(M,x)\otimes_{\Bbb F}\pi_j(M,x) \to\pi_{i+j-1}(M,x)$$
making  $\pi_{\bullet+1}(M,x)$ into a graded Lie algebra over 
$\Bbb F$. 
\endproclaim

This statement, see, e.g., \cite{Ka1, Prop. 1.2.2},
is in fact equivalent to the very definition of a homotopy Lie
algebra and should be regarded as being as old as this definition. 
More precisely \cite{St2}, if $\frak g^\bullet$ is a graded vector 
space, a structure of
a homotopy Lie algebra on $\frak g^\bullet$ is the same as a continuous
derivation $D$ of the completed symmetric algebra $\widehat S^\bullet(\frak 
g^*[-1])$ satisfying $D^2=0$ (so $(\widehat S^\bullet(\frak g^*[-1]), D)$
serves as the cochain complex of $\frak g^\bullet$). 
If we take $\frak g^\bullet = T^\bullet_xM[-1]$, then
the completed local ring $\widehat{\Cal O}_{M,x}^\bullet$ with its
natural differential $d$, serves as a cochain complex for $\frak g^\bullet$:
a choice of formal coordinates identifies it with $\widehat S^\bullet
(\frak g^*[-1])$.

These two Lie algebra structures (one on the tangent sheaf, the other
on its shifted fiber) will be related in (2.7.9).

As usual, any morphism $f: M\to N$ of dg-manifolds induces
a morphism of coherent dg-sheaves
$d^*f: f^*\Omega^{1\bullet}_N\to\Omega^{1\bullet}_M$ and
for any
$x\in M(\Bbb F)=\pi_0(M)(\Bbb F)$, a morphism of  complexes 
$d_xf: T_x^\bullet M\to T^\bullet_{f(x)} N$. These morphisms
of complexes fit together into a morphism of quasicoherent dg-sheaves
$df: T^\bullet M\to f^* T^\bullet N$.

The equivalence (i)$\Leftrightarrow$(ii) in the following 
proposition (see \cite{Ka1, Prop. 1.2.3}) 
can be seen as an analog of the Whitehead theorem
in topology.

\proclaim {(2.5.9) Proposition} (a) Let $f: M\to N$ be a morphism of 
dg-manifolds.
Then the following conditions are equivalent:

(i) $f$ is a quiasiisomorphism.

(ii) The morphism of schemes $\pi_0(f): \pi_0(M)\to \pi_0(N)$ is an 
isomorphism,
and for any field extension $\Bbb K\subset\Bbb F$
and any $\Bbb F$-point $x$ of $M$ the differential $d_xf$ induces
an isomorphism $\pi_i(M,x)\to \pi_i (N, f(x))$ for all $i< 0$. 

(iii) $\pi_0(f)$ is an isomorphism and $d^*f: f^*\Omega^{1\bullet}_N\to\Omega^{
1\bullet}_M$ is a quasiisomorphism of
coherent dg-sheaves on $M$.

(b) If any of these conditions is satisfied, then
$df: T^\bullet M\to f^* T^\bullet N$ is a quasiisomorphism. 

\endproclaim

\demo { Proof} We first establish the equivalences in (a)

(ii)$\Rightarrow$(i) To show that $f$ is an isomorphism,
it is enough to prove that for any field extension
$\Bbb F \supset \Bbb K$ and
any $x\in\pi_0(M)(\Bbb F)$
the map
$\hat f^*\widehat{\Cal  O}^\bullet_{N, f(x)}\to \widehat{\Cal  
O}^\bullet_{M,x}$
which $f$ induces on the completed local dg-algebras,
is a quasiisomorphism.
For that, notice that $\widehat{\Cal  O}^\bullet_{M,x}$ has a filtration
whose quotients are the symmetric powers of the cotangent dg-space
$T^*_xM$. So if $f$ gives a quasiisomorphism of tangent
dg-spaces, we find that $\hat f^*$ induces quasiisomorphisms on
the quotients of the natural filtrations. So the proof is
accomplished by invoking a spectral sequence argument, which
is legitimate (i.e., the spectral sequences converge) because
the dg-algebras in question are ${\Bbb Z}_{\leq 0}$-graded.

(iii)$\Rightarrow$(ii) Since $d^*f$ is a quasiisomorphism of dg-vector
bundles, it induces, by taking the tensor product with $\Bbb F_x$,
a quasiisomorphism on the fiber at each $x\in M(\Bbb K)$. The fibers
of $\Omega^{1\bullet}_M$ and $f^*\Omega^{1\bullet }_N$ at $x$
are just the complexes dual to $T^\bullet_xM$ and $T^\bullet_{f(x)}N$;
in particular, they are finite-dimensional in each degree.
Thus the dual map, which is $d_xf$, is a quasiisomorphism as well.

(i)$\Rightarrow$(iii) It is enough to show that the morphism
$d^\star f: f^{-1}\Omega^{1\bullet}_N\to\Omega^{1\bullet}_M$ is
a quasiisomorphism. Indeed, knowing this, (iii) is obtained by
applying the functor $-\otimes_{f^{-1}\Cal O_N^\bullet} 
\Cal O_M^\bullet$ (to pass from $f^{-1}$ to $f^*$) and invoking
the Eilenberg-Moore spectral sequence and the fact that
$f^{-1}\Cal O_N^\bullet \to\Cal O_M^\bullet$ is a quasiisomorphism. 
Since we can work locally, all we need to prove is a statement about
dg-algebras. We call a commutative dg-algebra smooth if its spectrum
is an affine dg-manifold. 

\enddemo

\proclaim {(2.5.10) Lemma} If $\phi: A\to B$ is a quasiisomorphism of
smooth $\Bbb Z_-$-graded commutative dg-algebras, then
$d\phi: \Omega^{1\bullet}_A\to \Omega^{1\bullet}_B$ is
a quasiisomorphism of complexes.
\endproclaim

\demo{Proof} This is a standard application of the theory of Harrison
homology, cf. \cite{Lo, \S 4.2.10}. The Harrison chain complex is
$$\text{Harr}_\bullet(A,A)  = \text{FCoLie}(A[-1])\otimes A$$
where FCoLie stands for the free graded coLie algebra generated by
a graded vector space\footnote{Since the primitive elements in a free
tensor coalgebra give the free coLie algebra, this description coincides
with the more traditional one which gives the Harrison complex
as the primitive elements in the Hochshild complex of $A$.}. It
satisfies the following properties:

\vskip .1cm

(a)  $\text{Harr}_\bullet(A,A)$ is covariantly functorial in $A$
and its dependence on $A$ is exact (takes quasiisomorphisms to
quasiisomorphisms).

(b) If $A$ is smooth and $d_A=0$, then $\text{Harr}_\bullet(A,A)$ is 
quasiisomorphic to
$\Omega^{1\bullet}_A$.

\vskip .1cm

Part (a) is obvious from the tensor nature of the functor FCoLie.
Part (b) can be proved in the same way as for ordinary (not dg)
smooth algebras: by realizing $\text{Harr}_\bullet(A,A)$ as
indecomposable elements in the Hochschild complex and using the
Hochschild-Kostant-Rosenberg theorem, see \cite{Lo, Thm. 3.4.4}. 
\enddemo

This concludes the proof of part (a) of Proposition 2.5.9. 
To prove (b), it is enough to work locally, so to
assume $M=\op{Spec}(A)$ is affine. Then we view $d^*f: f^*\Omega^{1\bullet}_N
\to \Omega^{1\bullet}_M$ as a morphism of dg-modules over $A$, and
we know that is a quasiisomorphism. 
 By further localizing on
$\op{Spec}(A^0)$, we can assume that $A^0$ is a local
ring. Then, projective modules over $A^0$ being free,
we have that both $f^*\Omega^{1\bullet}_N$ and $\Omega^{1\bullet}_M$
are quasifree dg-modules over $A$.  Hence, by Corollary (2.1.5),
the dual morphism to $d^*f$, i.e., $df: T^\bullet M\to f^*T^\bullet N$
is a quasiisomorphism as well.

\vskip .3cm

\noindent {\bf (2.6) Existence of smooth resolutions.}
A morphism of dg-schemes $f: X\to Y$ will be called a closed
embedding, if $f_0: X^0\to Y^0$ is a closed embedding of
schemes and the structure morphism of sheaves of dg-algebras
$f_0^*\Cal  O^\bullet_N\to\Cal  O^\bullet_M$
is surjective. 

\proclaim{(2.6.1) Theorem} (a) For any quasiprojective dg-scheme $X$ there
is a dg-manifold $M$ and a quasiisomorphic closed embedding $X\hookrightarrow
M$.

(b) Given any two embeddings $X\hookrightarrow M$, $X\hookrightarrow N$
as in (a), they can be complemented by quasiisomorphic closed embeddings
$M\hookrightarrow L$, $N\hookrightarrow L$ for some dg-manifold $L$ so
that the natural square is commutative.

\endproclaim

\demo{Proof} (a)$\Rightarrow$(b). Given $M,N$, we set
$$Y=M\cup_XN = \biggl( M^0\cup_{X^0} N^0, \,\, \Cal  O_M^\bullet\times
_{\Cal  O_X^\bullet} \Cal  O_N^\bullet\biggr).$$
Then we have a diagram as required except that $Y$ may be not
smooth. To amend this, it suffices to embed $Y$ into a dg-manifold $L$
as in (a). 

\vskip .1cm

(a) This is a version of the standard fact asserting the existence of a free
resolution for a $\Bbb Z_-$-graded dg-algebra, see, e.g., \cite{BG}
\S 4.7.  If $X$ is affine, then this fact indeed implies
what we need. 

\enddemo

In the general case, let $X=(X^0, \Cal  O_X^\bullet)$ be given.
As $X^0$ is a quasiprojective scheme, we can choose its embedding into a $P^n$
as a locally closed subscheme. Take $M^0$ to be an open subset in $P^n$
such that $X^0$ is closed in $M^0$. We then construct $\Cal  O_M^\bullet$
by induction as the union of sheaves of dg-subalgebras
$$\Cal  O_{M^0} = \Cal  O_0^\bullet \subset \Cal  O_1^\bullet \i\Cal  
O_2^\bullet ...$$
such that:

\vskip .1cm

\item{(1)} $\Cal  O_i^\bullet$ is obtained from $\Cal  O_{i-1}^\bullet$
 by adding a vector bundle of generators in degree $(-i)$.

\vskip .1cm

\item{(2)} We have a compatible system of algebra morphisms 
$p^{(i)}: \Cal  O_i^\bullet \to\Cal  O_X^\bullet$ so that each $p^{(i)}$
is bijective on $\underline{H}^j$ for $-i+1\leq j\leq 0$
and surjective on the sheaf of $j$th cocycles for $-i\leq j\leq 0$. 

\vskip .1cm

The following elementary lemma shows that surjectivity on cocycles
implies surjectivity on graded components and therefore the
map we will construct in this way will be a closed
embedding. 

\proclaim {(2.6.2) Lemma} Let $\phi: C^\bullet\to D^\bullet$
be a quasiisomorphism of complexes of vector spaces which is
surjective on cocycles. Then each $\phi_i: C^i\to D^i$ is surjective.
\endproclaim

The inductive construction of the $\Cal O^\bullet_i$ follows 
the standard pattern of
``imitating the procedure of attaching cells to kill homotopy groups"
(\cite {Q2}, p. 256, see also \cite{BG}, \S 4.7)
  except that we use Lemma 2.3.4 to produce
a vector bundle of generators. We leave the details to the reader.

\vskip .3cm

\noindent {\bf (2.7) Smooth morphisms.} We now relativize the above
discussion.

\proclaim{(2.7.1) Definition} Let $M, N$ be dg-schemes of finite type.
A morphism $f: M\to N$ is called smooth, if the following
two conditions hold:

(a) The underlying morphism $f_0: M^0\to N^0$ of ordinary schemes is smooth.

(b) Locally, on the Zariski topology of $M^0$, we have an isomorphism
of graded algebras
$$\Cal O^\bullet_{M\#} \simeq S_{M^0}(Q^\bullet)\otimes f_0^* \Cal O^\bullet_{
N\#},$$
where $Q^\bullet = \bigoplus_{i\leq -1}Q^i$ is a graded vector
bundle on $M^0$.

\endproclaim

As before, for a smooth morphism we can always globally define the graded
 bundle
$$Q^\bullet_{M/N} \quad = \biggl( \Cal O_{M\#}^{\leq -1} \bigl/
(\Cal O_{M\#}^{\leq -1})^2\biggr) \otimes_{f_0^*\Cal O_{N\#}^\bullet}
\Cal O_M^\bullet,\leqno (2.7.2)$$
but it embeds into $\Cal O_M^\bullet$ only locally. We also have the relative
cotangent dg-sheaf
$\Omega^{1\bullet}_{M/N}$ which is a $\Bbb Z_-$-graded vector bundle
and the relative tangent dg-sheaf
$$T^\bullet(M/N) = \Cal Der^\bullet_{f_0^{-1}\Cal O_N^\bullet}
(\Cal O_M^\bullet, \Cal O_M^\bullet) = \Cal Hom _{\Cal O_M^\bullet}
(\Omega^{1\bullet}_{M/N}, \Cal O_M^\bullet).\leqno (2.7.3)$$
This is a quasicoherent dg-sheaf.

Let $x: N\to M$ be a section of $f$ (i.e., an $N$-point of an $N$-dg-scheme
$M$). Then we define the tangent dg-space (or bundle) to $M/N$ at
(or along) $x$ as
$$T^\bullet_x(N/M) = x^* T^\bullet (N/M) = x_0^* T^\bullet(N/M)
\otimes_{x^*_0 \Cal O_M^\bullet} \Cal O_N^\bullet.\leqno (2.7.4)$$
This is a quasicoherent dg-sheaf on $N$. 

\vskip .1cm

\noindent {\bf (2.7.5) Remarks.} (a) As in (2.5.7), one can show that 
$T^\bullet_x(M/N)[-1]$ is a ``sheaf of homotopy Lie algebras" on $N$;
in particular, its cohomology $\underline{H}^\bullet T^\bullet_x(M/N)[-1]$
is naturally a sheaf of graded Lie algebras. Indeed, the role  of the
``cochain complex" of $T^\bullet_x(M/N)[-1]$  is played by 
$\widehat{\Cal O}^\bullet_{M, x}$
the completion of $\Cal O_M^\bullet$ along the subscheme $x(N)$.

(b) Globally it may be impossible to identify 
$\widehat{\Cal O}^\bullet_{M, x}$ with the symmetric algebra of 
 $T^\bullet_x(M/N)[-1]$ and the corresponding obstruction gives rise
to another Lie algebra-type structure,
present even when $M,N$ are ordinary (not dg) schemes. More precisely, the 
obstruction
to splitting the second infinitesimal neighborhood gives rise to a version 
of the Atiyah class:
$$\alpha\in H^1(N^0, \Cal Hom (S^2 T^\bullet_x(M/N)[-1],  T^\bullet_x(M/N)[-1])
$$
 which formally satisfies the Jacobi identity, as an element of
an appropriate operad. This generalizes the
main observation of   \cite{Ka2, Thm 3.5.1}, which
corresponds to the case when  $N=X$ is an ordinary manifold, $M=X\times X$
and $x$ is the diagonal map. 

\vskip .1cm

The following smoothing statement can be regarded as a rudiment
of a closed model structure in the category of dg-schemes.

\proclaim{(2.7.6) Theorem} Let $f: M\to N$ be any morphism
of quasiprojective dg-schemes. Then $f$ can be factored
as $M\buildrel i\over\hookrightarrow \widetilde M
\buildrel \widetilde f\over\rightarrow N$, where $i$ is a quasiisomorphic
closed embedding and $\widetilde f$ is smooth. Any two such factorizations can 
be
included into a third one.
\endproclaim

\demo{Proof} We embed $M_0$ into $N_0\times P^n$.
Then take for $\widetilde M_0$ an open subset in $N_0\times P^n$
such that $M_0$ is closed in $\widetilde M_0$.
Then the procedure is the same as outlined in (2.6.1) for $N=pt$.
\enddemo
A diagram $M\buildrel i\over\hookrightarrow \widetilde M
\buildrel \widetilde f\over\rightarrow N$ as in (2.7.6) will be
called a smooth resolution of $f: M\to N$.

The following is a relative version of a part of Proposition 2.5.9, proved
in the same way. 

\proclaim{(2.7.7) Proposition} Let 
$$\matrix M_1&\buildrel q\over\longrightarrow &M_2\\
f_1\searrow&&\swarrow f_2\\
&N&\\
\endmatrix$$
be a commutative triangle with $f_i, i=1,2$ smooth and $q$ a quasiisomorphism.
Then $dq: T^\bullet(M_1/N)\to q^*T^\bullet(M_2/N)$ and $d^*q:
q^*\Omega^{1\bullet}_{M_2/N}\to\Omega^{1\bullet}_{M_1/N}$ are
quasiisomorphisms. 
\endproclaim

\proclaim {(2.7.8) Definition} The derived relative tangent complex of
a morphism $f: M\to N$ of quasiprojective dg-schemes is, by definition,
$$RT^\bullet(M/N) = T^\bullet(\widetilde M/N),$$
where $\widetilde f:\widetilde M\to N$ is any smooth resolution of $f$.
This is a sheaf of dg-Lie algebras on $\widetilde M$. 
\endproclaim

Proposition 2.7.7 guarantees that $RT^\bullet(M/N)$ is well defined up to
quasiisomorphism. 

\vskip .1cm

\noindent {\bf (2.7.9) Example.} Let $M=\{x\}$ be a $\Bbb K$-point of
$N$ and $f$ be the embedding of this point.  Then the derived
relative tangent complex $RT^\bullet(\{x\}/N)$  is quasiisomorphic
to the shifted tangent complex $T^\bullet_xN[-1]$. This can be
seen by taking for $\widetilde M$ the spectrum of 
a Koszul resolution of $\Bbb K_x$ on an affine open
dg-submanifold $U\i N$ containing $x$. 
The presence of the Lie bracket on $T^\bullet(\widetilde M/N) = RT^\bullet
(\{x\}/N)$
  provides an alternative explanation of the
presence of a  homotopy Lie algebra structure on $T^\bullet_xN[-1]$.
 
\vskip .3cm

\noindent {\bf (2.8) Derived fiber products.} Let $f_i: M_i\to N$
be morphisms of dg-schemes, $i=1,2$. The fiber product
$M_1\times_N M_2$ is defined as follows. First, we  form
the fiber product of underlying ordinary schemes:
$$\minCDarrowwidth{7 mm}
\CD
 M_1^0\times_{N^0} M_2^0@>g_{2,0}>> M_2^0\\
@Vg_{1,0}VV @VV f_{2,0} V\\
M_1^0 @>>f_{1,0}> 
N^0
\endCD
$$
and define 
$$M_1\times_N M_2\quad =\quad \biggl( M_1^0\times_{N^0} M_2^0, \,\,
g_{2,0}^{-1} \Cal O_{M_2}^\bullet
 \otimes_{(f_{1,0}g_{1,0})^{-1}\Cal O_N^\bullet}
g_{1,0}^{-1}\Cal O_{M_1}^\bullet\biggr),$$
so that we have the natural square
$$\minCDarrowwidth{7 mm}
\CD
 M_1\times_{N} M_2@>g_{2}>> M_2\\
@Vg_{1}VV @VV f_{2} V\\
M_1 @>>f_{1}> N
\endCD\leqno (2.8.1)
$$

The following fact is clear.

\proclaim{(2.8.2) Proposition} If $f_2$ is a smooth morphism, then
so is $g_1$.
\endproclaim

The fiber (or preimage) is a particular case of this construction.
More precisely, let $f: M\to N$ be a morphism and $y\in N(\Bbb K)$ be
a point. The fiber $f^{-1}(y)$ is the fiber product
$M\times_N \{y\}$. If $f$ is a smooth morphism, then $f^{-1}(y)$ is
a dg-manifold. Suppose further that $M,N$ and $f$ are all smooth. 
Then we have the Kodaira-Spencer map
$$\varkappa: T^\bullet_yN[-1]\to R\Gamma(f^{-1}(y), T^\bullet f^{-1}(y))
\leqno (2.8.3)$$
which is, as in the standard case, obtained from the short exact sequence
$$0\to T^\bullet(M/N)\to T^\bullet M\to f^*T^\bullet N\to 0$$
by tensoring over $\Cal O_M^\bullet$ with $\Cal O_{f^{-1}(y)}^\bullet$
and using the adjunction. 

\vskip .1cm

\noindent {\bf (2.8.4) Remark.} Note that both the source and target
of $\varkappa$ possess a homotopy Lie algebra structure: the source by
Proposition 2.5.8, the target as the direct image of a sheaf of dg-Lie
 algebras.
In fact, it can be shown that $\varkappa$ is naturally a homotopy morphism
of homotopy Lie algebras. In particular, the graded Lie algebra
$\pi_{\bullet+1}(N, y)$ acts on the hypercohomology  space
 $\Bbb H^\bullet(f^{-1}(y), \Cal O^\bullet)$ in a way remindful of the 
monodromy action of a fundamental
group. We postpone further discussion to a more detailed exposition
of the basics of the theory, to be completed at a future date. 

\vskip .1cm

Note that a smooth morphism is flat (this is proved in the same way
as for the case of usual schemes). Therefore Proposition 2.4.1
implies the following.

\proclaim{(2.8.5) Proposition} Suppose that in the fiber product diagram
(2.8.1) the morphism $f_1$ (and hence $g_2$) is smooth. Then we have
a converging 
(Eilenberg-Moore) spectral sequence of quasicoherent sheaves of graded 
algebras
on $M^0$:
$$E_2 = \Cal Tor_\bullet^{(f_{10}g_{10})^{-1}
\underline{H}^\bullet(\Cal O^\bullet_N)}
\bigl( g_1^{-1}\underline{H}^\bullet(\Cal O^\bullet_{M_1}), \,\,
g_2^{-1}\underline{H}^\bullet(\Cal O_{M_2})\bigr) \quad
\Rightarrow \quad \underline{H}^\bullet(\Cal O_{M_1\times_N M_2}^\bullet).$$
\endproclaim

\proclaim {(2.8.6) Definition} The derived (or homotopy) fiber product
$M_1\times_N^R M_2$ is defined as $\widetilde M_1\times_N M_2$ where
$\widetilde f_1: \widetilde M_1\to N$ is a smooth resolution of 
the morphism $f_1$. 
\endproclaim

\proclaim{(2.8.7) Lemma} The definition of $M_1\times_N^R M_2$
is independent, up to quasiisomorphism, of the choice of a smooth
resolution of $f$.

\endproclaim

\demo{Proof} It is enough to show that whenever we have a diagram
$$\matrix \widetilde M_1& \buildrel q\over \hookrightarrow& \widetilde M'_1\\
\widetilde f_1 \searrow&&\swarrow\widetilde f_1'\\
&N&\\
\endmatrix$$
with $q$ a quasiisomorphism and $\widetilde f_1, \widetilde f'_1$ smooth,
the induced morphism $\widetilde M_1\times_N M_2\to \widetilde M_1'\times_N
M_2$ is a quasiisomorphism. For this notice that we have a morphism from
the Eilenberg-Moore spectral sequence calculating  $\underline{H}^\bullet
\bigl(\Cal O^\bullet_{\widetilde M_1\times_N M_2}\bigr)$
to the similar sequence calculating 
$\underline{H}^\bullet
\bigl(\Cal O^\bullet_{\widetilde M'_1\times_N M_2}\bigr)$ and this
morphism is an isomorphism on $E_2$ terms. 
\enddemo 

We can now formulate the final form of the Eilenberg-Moore spectral
 sequences for the
derived fiber products.

\proclaim {(2.8.8) Proposition} Suppose a square of quasiprojective
dg-schemes
$$\minCDarrowwidth{7 mm}
\CD
 M @>g_{2}>> M_2\\
@Vg_{1}VV @VV f_{2} V\\
M_1 @>>f_{1}> N
\endCD
$$
is homotopy cartesian, i.e., the natural morphism $m\to M_1\times_N^R M_2$
is a quasiisomorphism. Then we have the two convergent  spectral sequences
 of sheaves of algebras  on $M^0$:
$$E_1 = \Cal Tor_\bullet^{(f_{1,0}g_{1,0})^{-1}(\Cal O^\bullet_{N\#})}
\bigl( g_1^{-1}\Cal O^\bullet_{M_1\#},\,\, g_2^{-1}\Cal O^\bullet_{M_2\#}
\bigr) \quad\Rightarrow\quad \underline {H}^\bullet(\Cal O_M^\bullet),$$
$$E_2  =  \Cal Tor_\bullet^{(f_{1,0}g_{1,0})^{-1}(\underline{H}^\bullet (\Cal 
O^\bullet_{N}))}\bigl( g_1^{-1}\underline{H}^\bullet(\Cal O^\bullet_{M_1}),
\,\, g_2^{-1} \underline{H}^\bullet(\Cal O^\bullet_{M_2})\bigr)
\quad\Rightarrow\quad
 \underline {H}^\bullet(\Cal O_M^\bullet).$$

\endproclaim

\noindent {\bf (2.8.9) Remark.} More generally, one can
define the derived fiber product for any morphisms $f_i: M_i\to N$
of arbitrary dg-schemes (not necessarily quasiprojective or of finite type).
But we need to assume that at least one of the $f_i$ can be quasiisomorphically
replaced by a $\#$-flat morphism $\widetilde M_i\to N$. This is
the case, for example, when $f_i$ is an affine morphism 
(use the relative bar-resolution).

\vskip .1cm

\noindent {\bf (2.8.10) Examples} {\it (a) Derived intersection.} 
If $Y,Z$ are closed subschemes of a quasiprojective dg-scheme  $X$, then
the derived intersection $Y\cap^R Z$ is defined as the derived
fiber product of $Y$ and $Z$ over $X$. If $X,Y,Z$ are ordinary
(not dg) schemes, then the cohomology sheaves of $\Cal O^\bullet_{Y\cap^R Z}$
are the $\Cal Tor_i^{\Cal O_X}(\Cal O_Y, \Cal O_Z)$, see \cite{Kon}, n.
(1.4.2). 

\vskip .1cm

{\it (b) Homotopy fibers}. Given any morphism $f: M\to N$ of quasiprojective
dg-schemes and any point $y\in N(\Bbb K)$, we have the homotopy fiber
$Rf^{-1}(y):= \widetilde f^{-1}(y)$ where $\widetilde f$ is a smooth
resolution of $f$. Note that for $N$ smooth (and $f$ arbitrary)
we always have the derived Kodaira-Spencer map
$$R\varkappa: T^\bullet_y N[-1]\to R\Gamma(Rf^{-1}(y), T^\bullet)$$
and Remark 2.8.4 applies to this situation as well.  

\vskip .1cm

{\it (c) The loop space.} Consider the particular case of (b), namely
$M=\{y\}$ and $f = i_y$ being the embedding. Using the topological analogy,
it is natural to call the homotopy fiber $Ri_y^{-1}(y)$ the
loop space of $N$ at $y$ and denote it $\Omega(N,y)$. This dg-scheme has
only one $\Bbb K$-point, still denoted $y$ (``the constant loop"). As
for the tangent space at this point, 
we have $T^\bullet_y \Omega(N,y) = T^\bullet_yN [-1]$
and the derived Kodaira-Spencer map for $i_y$ is the identity. 

By going slightly beyond the framework of this paper, we can
make the analogy with the usual loop space even more pronounced.
Namely, consider the $\Bbb Z_+$-graded dg-algebra $\Lambda[\xi]$,
$\text{deg}(\xi)=1$, in other words, $\Lambda[\xi] = H^\bullet(S^1, 
\Bbb K)$ is the topological cohomology of the usual circle. Let us
formally associate to this algebra the dg-scheme $\Bbb S = 
\text{Spec}(\Lambda[\xi])$ (``dg-circle"). It has a unique $\Bbb K$-point
which we denote $e$. Then we can identify $\Omega(N,y)$ with the
internal Hom in the category of pointed dg-schemes
$$\Omega(N,y) = \underline{\text{Hom}}\bigl( (\Bbb S, e), \, (N, y)\bigr),$$
similarly to the usual definition of the loop space. 

Further, the fact that the usual loop space is a group up to homotopy, has
the following analog, cf. \cite{Q1}. 
Let $\Pi\to N$ be a smooth quasiisomorphic
replacement of $i_y: \{y\}\to N$, see (2.7.9). Then we have a groupoid
$\frak G$ in the category of dg-schemes with
$$\text{Ob}\,\,\frak G = \Pi \sim \{\text{pt}\}, \quad \text{Mor}\,\,
\frak G = \Pi\times_N\Pi \sim \Omega(N, y).$$
This group-like structure on $\Omega(N,y)$ provides still
another explanation of the fact that its tangent space
$T^\bullet_y N[-1]$ is a homotopy Lie algebra.

\vfill\eject

\heading { 3. A finite-dimensional model}
\endheading

\subhead{(3.1) The problem}\endsubhead 
Our goal in this paper is to construct, in the situation of (1.1),
a dg-manifold $RSub_h(\Cal  F)$ satisfying the conditions (0.3.1) and
(0.3.2). In this section we consider a finite-dimensional analog of this
problem. Namely, let $A$ be a finite-dimensional associative
algebra, $M$ a finite-dimensional
left $A$-module and $G_A(k, M)$ the $A$-Grassmannian, see (1.3). We want
  to construct a dg-manifold $RG_A(k,M)$
with the properties;
$$\pi_0 \, RG_A(k,M) = G_A(k,M), \quad H^i T^\bullet_{[V]}RG_A(k,M) = 
\operatorname{ Ext}^i_A(V, M/V).\leqno (3.1.1)$$
 As we will see later, the problem of constructing the derived 
Quot scheme can  be reduced to this.

\vskip .3cm

\subhead{(3.2) Idea of construction}\endsubhead We first 
realize $G_A(k,M)$ inside
the (noncompact) smooth variety 
$|\Cal Hom(A\otimes \widetilde V,
\widetilde V)|$ as follows. The bundle 
$\Cal Hom(A\otimes \widetilde V,
\widetilde V)$ has a canonical section $s$ defined over the subscheme
$G_A(k,M)$. This section is given by the $A$-action $A\otimes V\to V$
present on any submodule $V\subset M$. We embed $G_A(k,M)$
into  $|\Cal Hom(A\otimes \widetilde V,
\widetilde V)|$ as the graph of this section and will construct $RG_A(k,M)$
so that its underlying ordinary manifold is 
$|\Cal Hom(A\otimes \widetilde V,
\widetilde V)|$. For this, we will represent the embedded $G_A(k,M)$
as the result of applying the following two abstract constructions.

\vskip .2cm

\subhead{(3.2.1) The space of actions}\endsubhead 
Let $V$ be a finite-dimensional
vector space. Then we have the subscheme $\operatorname{ Act}(A,V)$ in the
affine space $\operatorname{ Hom}_\Bbb K(A\otimes_\Bbb K V, V)$ consisting
of all $A$-actions (i.e., all $A$-module structures) on $V$. Note
that we do not identify here two $A$-module structures which
give isomorphic modules.

\vskip .2cm

\subhead{ (3.2.2) The linearity locus}\endsubhead Let $S$ be a scheme, 
and $M,N$ be two
vector bundles over $S$ with $A$-actions in fibers. In other words,
$M,N$ are ${\Cal  O}_S\otimes_\Bbb KA$-modules which are locally
free as ${\Cal  O}_S$-modules. Let also $f: M\to N$ be an ${\Cal  O}_S$-linear
morphism. Its linearity locus is the subscheme
$$\operatorname{ Lin}_A(f) = \biggl\{ s\in S\bigl| f_s: M_s\to N_s\quad 
\operatorname{ is}
\quad A\operatorname{ -linear}\biggr\}.$$
This is just the fiber product
$$\minCDarrowwidth{7 mm}
\CD
\operatorname{ Lin}_A(f)@>>> S\\
@VVV @VV f V\\
|\Cal Hom_{A\otimes{\Cal   O}_S}(M,N)| @>>> 
|\Cal Hom_{{\Cal   O}_S}(M,N)|
\endCD$$

\vskip .2cm

Let us apply the first construction to each fiber $V$ of the bundle
$\widetilde V$ on $G(k,M)$. We get the fibration
$$\operatorname{ Act}(A, \widetilde V)\buildrel q\over\longrightarrow G(k,M)$$
which is embedded into $|\Cal Hom(A\otimes\widetilde V, \widetilde V)|$.
By construction, the pullback $q^* \widetilde V$ is a bundle of $A$-modules.
Let also $\underline{M}$ be the  trivial bundle of $A$-modules on
$\operatorname{ Act}(A, \widetilde V)$ with fiber $M$. Then, we have the 
tautological
morphism $f: q^*\widetilde V\to\underline{M}$ of vector bundles whose
fiber over a point $([V], \alpha)\in \operatorname{ Act}(A, \widetilde V)$ is 
just
the embedding $V\hookrightarrow M$.

\proclaim {(3.2.3) Proposition} For this morphism $f$ the scheme
$\operatorname{ Lin}_A(f)\subset \operatorname{ Act}(A,\widetilde V)$ 
coincides with
$G_A(k,M)$ embedded into
$\operatorname{ Act}(A, \widetilde V)\subset 
|\Cal Hom(A\otimes \widetilde V,
\widetilde V)|$.\endproclaim

\demo{Proof} Given a linear subspace $V\subset M$ and an $A$-action
$\alpha: A\otimes V\to V$, the condition that the embedding $V\hookrightarrow
M$ be $A$-linear precisely means that $V$ is a submodule and $\alpha$
is the induced action.\enddemo

Now the idea of constructing $RG_A(k,M)$ 
is to develop the derived analogs of
the two constructions (3.2.1), (3.2.2) and apply them to the situation
just described.

\vskip .3cm

\subhead{(3.3) The derived space of actions}\endsubhead Let us first analyze 
in more detail the
(non-derived) construction $\operatorname{ Act}(A,V)$. It can be defined
for any (possibly infinite-dimensional) associative algebra $A$ and
a finite-dimensional vector space $V$. In this case we can
apply the conventions of (2.2.3)(b) to the complex 
$F=\operatorname{ Hom}(V, A\otimes V)$
and denote 
the affine scheme $|F^*| = \operatorname{ Spec}\, S(F)$ by 
$|\operatorname{ Hom}(A\otimes V, V)|$.
The scheme  $\operatorname{ Act}(A,V)$ is 
the closed subscheme of $|\operatorname{ Hom}(A\otimes V, V)|$, 
whose coordinate
ring is $S(\operatorname{ Hom}(V, A\otimes V))$ modulo the ideal 
expressing the
associativity conditions. At the level of $\Bbb K$-points,
$\mu: A\otimes V\to V$ is an
action if and only if the map
$$\delta\mu: A\otimes A\otimes V\to V,\quad a_1\otimes a_2\otimes v\mapsto
\mu(a_1a_2\otimes v)-\mu(a_1\otimes\mu(a_2\otimes v)),\leqno (3.3.1)$$
vanishes.  In this case $T_\mu \operatorname{ Act}(A,V)$ is identified
with the space of 1-cocycles in the bar-complex
$$\operatorname{ Hom}_\Bbb K(V,V)\to \operatorname{ Hom}_\Bbb K(A\otimes V, V)
\to \operatorname{ Hom}_\Bbb K(A\otimes A\otimes V, V)\to ...\leqno (3.3.2)$$
calculating $\operatorname{ Ext}^\bullet_A(V, V)$.

\vskip .5cm

\noindent {\bf (3.3.3) Remark.} The reason that we get the
space of 1-cocyles instead of the cohomology which is a more invariant
object is that we do not identify isomorphic module structures.
If we consider the quotient stack of 
$\operatorname{ Act}(A,V)$ by $GL(V)$, then
the tangent space to this stack at a point $\mu$ is a 2-term
complex concentrated in degrees $0, -1$ and
$$H^i T^\bullet_\mu\biggl( \operatorname{ Act}(A,V)/GL(V)\biggr) 
= \operatorname{ Ext}^{i+1}_A(V,V), \quad i=0, -1.$$

Our aim in this subsection is to construct, for each finite-dimensional
$A$,  a (smooth) dg-manifold 
$\operatorname{ RAct}(A,V)$ with $\pi_0 = \operatorname{ Act}(A,V)$ and the
tangent space at any $\mu\in \operatorname{ Act}(A,V)$ having
$$H^iT^\bullet_\mu \operatorname{ RAct}(A,V) = 
\cases T_\mu\operatorname{ Act}(A,V),& i=0\\
\operatorname{ Ext}^{i+1}_A(V,V),&  i>0\endcases
\leqno (3.3.4)$$

The method of construction will be the standard approach of
homological algebra, namely using free associative resolutions of $A$.
This is similar to C. Rezk's approach \cite{Re} to constructing
``homotopy" moduli spaces for actions of an operad.  More precisely,
we will construct for any, possibly
infinite-dimensional $A$, an affine dg-scheme
$\operatorname{ RAct}(A,V)$ whose coordinate algebra is free commutative,
and will show that for $\dim(A)<\infty$, we can choose a representative
with finitely many free generators in each degree, so that
we have a dg-manfold.

 \vskip .2cm

Notice first that the  construction of $\operatorname{ Act}(A,V)$ in the 
beginning
of this subsection can
be carried through for any 
  $\Bbb Z_-$-graded associative dg-algebra
(with $A^i$ possibly infinite-dimensional) and $V$ 
a finite-dimensional vector space 
(which we think as being graded, of degree 0). As in the ungraded
case, $\operatorname{ Act}(A,V)$ is a closed dg-subscheme in 
$|\operatorname{ Hom}(A\otimes V, V)|$ given by the ideal
of associativity conditions, which is now a dg-ideal. 
The association $A\mapsto \operatorname{ Act}(A,V)$ is functorial:
a morphism of dg-algebras $f:A_1\to A_2$ gives rise to a morphism
of dg-schemes $f^*: \operatorname{ Act}(A_2, V)\to \operatorname{ Act}(A_1, 
V)$.

Next, assume that $A=F(E^\bullet)$ is a free associative (tensor)
algebra without unit generated by a $\Bbb Z_-$-graded vector space 
$E^\bullet$. Then,
clearly, we have
$$\operatorname{ Act}(F(E^\bullet), V) = |\operatorname{ Hom}_\Bbb 
C(E^\bullet\otimes V, V)|,\leqno (3.3.5)$$
as an action is uniquely defined by the action of generators,
which can be arbitrary.

 Further, assume that $B$ is a $\Bbb Z_-$-graded associative dg-algebra
which is quasifree, i.e., 
such that $B_\#  \simeq F(E^\bullet)$ is free.
Then, the graded scheme $\operatorname{ Act}(B,V)_\# $ is, by the above,
identified with $|\operatorname{ Hom}_\Bbb K(E\otimes V, V)|$.

\proclaim {(3.3.6) Proposition} If $f: B_1\to B_2$ is a quasiisomorphism
of quasifree associative $\Bbb Z_-$-graded dg-algebras, then
 $f^*: \operatorname{ Act}(B_2, V)\to \operatorname{ Act}(B_1, V)$ 
is a quasiisomorphism of
dg-schemes.\endproclaim

We will prove this proposition a little later. Assuming it is true,
we give the following definition.

\proclaim {(3.3.7) Definition} We define 
$\operatorname{ RAct}(A,V) = \operatorname{ Act}(B,V)$,
where $B\to A$ is any quasifree resolution.\endproclaim

\vskip .3cm

\subhead{(3.4) Reminder on $A_\infty$-structures}\endsubhead In what follows
it will be convenient to use the language of $A_\infty$-structures.
This concept goes back to J. Stasheff \cite{St1} for dg-algebras,
but here we need the companion concepts for modules 
(introduced by M. Markl \cite{Ma}) and
for morphisms of modules.

\proclaim {(3.4.1) Definition}
Let $A$ be an associative dg-algebra. A left $A_\infty$-module
over $A$ is a graded vector space $M$ together with $\Bbb K$-multilinear maps
$$\mu_n: A^{\otimes n}\otimes M\to M,\quad n\geq 0,
\quad \operatorname{ deg}(\mu_n)=1-n,$$
satisfying the conditions:
$$\sum_{i=1}^n (-1)^{\bar a_1 +...+\bar a_{i-1}} \mu_n(a_1,
..., da_i, ..., a_n, m) =$$
$$= \sum_{i=1}^{n-1} (-1)^i \mu_{n-1}(a_1, ..., a_ia_{i+1},
..., a_n,  m) -$$
$$-\sum_{p,q\geq 0\atop p+q=n} (-1)^{q(\bar a_1+...+\bar a_p) + 
p(q-1) + (p-1)q}
\mu_p(a_1,  ...,  a_p,  \mu_q(a_{p+1},  ..., a_n,  m)).$$\endproclaim

This implies, in particular, that $d_M=\mu_0$ satisfies $d_M^2=0$
and $\mu_1$ induces on $H^\bullet_{d_M}(M)$ a structure of a graded left
$H^\bullet (A)$-module.
 A collection of maps $\mu_n$ satisfying  the conditions of 
(3.4.1)
will be also referred to as an $A_\infty$-action of $A$ on $M$.
An $A_\infty$-action with $\mu_n=0$ for $n\geq 2$ is the same as
a structure of a dg-module in the ordinary sense.

\proclaim {(3.4.2) Definition} Let $A$ be an associative dg-algebra,
$M$ be a left $A_\infty$-module and $N$ be a genuine dg-module
over $A$. An $A_\infty$-morphism $f: M\to N$ is a collection
of linear maps
$$f_n: A^{\otimes n}\otimes M\to N, \quad \operatorname{ deg}(f_n) = -n,$$
satisfying the conditions:
$$df_n(a_0, ..., a_n, m)-\sum_{i=1}^n (-1)^i 
f_n(a_1, ..., da_i, ...,  a_n, m) =$$
$$\sum_{i=0}^{n-1} (-1)^i f_{n-1}(a_1, ..., a_ia_{i+1},
 ..., a_n, m) + \sum_{p=0}^n (-1)^{p(n-p)}
f_p(a_0, ..., a_p,  \mu_p(a_{p+1}, ...,
a_n, m)).$$\endproclaim

Again, the conditions imply that $f_0:M\to N$  is a morphism
of complexes and induces a morphism of left
$H^\bullet (A)$-modules $H^\bullet (M)\to H^\bullet (N)$.

\vskip .2cm

$A_\infty$-structures have transparent interpretation via bar-resolutions.
Let us start with $A_\infty$-modules. 
Assume that $A$ is $\Bbb Z_-$-graded and consider the graded vector space
 $$\bigoplus_{n=1}^\infty A^{\otimes n}[n-1]=
\operatorname{ Tot} \biggl\{ ... A\otimes A\otimes A ; \quad A\otimes A; \quad 
A\biggr\},
\leqno (3.4.3)$$
Here Tot means the $\Bbb Z_-$-graded vector space associated to a 
$\Bbb Z_-\times \Bbb Z_-$-graded
one. Let $D(A)$ be the free associative
algebra   without unit on this graded vector space.
The multiplication operation in  $D(A)$ will be denoted by $*$. We 
introduce a differential $d=d'+d''$ on $D(A)$  where $d'$ comes from
the tensor product differential on the $A^{\otimes m}$ and $d''$ is defined
on generators by
 $$d''(a_0\otimes ... \otimes a_n) = \sum_{i=0}^{n-1} (-1)^i a_0\otimes ... 
\otimes
a_ia_{i+1}\otimes ... \otimes a_n \ - \
\sum_{i=0}^{n-1} (-1)^i (a_0\otimes ... \otimes a_i)*(a_{i+1}\otimes ...
\otimes a_n) \leqno (3.4.4)$$

\proclaim {(3.4.5) Proposition} (a) The differential $d$ satisfies $d^2=0$.

\noindent (b) The projection $D(A)\to F(A)\buildrel m\over\rightarrow A$,
where $m$ is the multiplication in $A$, is a quasiisomorphism.\endproclaim

\demo{Proof} Well known: $D(A)$ is the bar-construction of
the cobar-construction of $A$, see \cite{HMS, II\S 3} for $(a)$
and \cite{HMS, Thm. II.4.4 } for $(b)$.\enddemo

\vskip .1cm

Thus $D(A)$ is a quasifree resolution of $A$.
By comparing (3.4.4) with (3.4.1), we find at once (cf. \cite{Ma}).

\proclaim{(3.4.6) Proposition} An $A_\infty$-action of $A$ on $M$ is the same
as a genuine action (structure of a dg-module) of $D(A)$ on $M$.\endproclaim

\vskip .1cm

Similarly, let $M$ be a left  $A_\infty$-module over $A$. We consider
the graded vector space 
$$ \operatorname{ Bar}_A(M) = \bigoplus_{n=1}^\infty A^{\otimes n}\otimes M 
[n-1] = \operatorname{ Tot}
 \biggl\{ ...\to A\otimes_\Bbb K A\otimes_\Bbb K\otimes M
\to A\otimes_\Bbb K M\biggr\}, \leqno (3.4.7)$$
cf. (2.4.4). 
It has a natural structure of a free left $A_\# $-module.
We equip it
with the differential
$$d(a_0\otimes ... \otimes a_n\otimes m) =
\sum_{i=0}^n (-1)^{i-1} a_0\otimes ... \otimes da_i\otimes ... \otimes a_n
\otimes m +
\leqno (3.4.8)$$
$$+\sum_{i=0}^{n-1} (-1)^i a_0\otimes ... \otimes a_ia_{i+1}\otimes
...\otimes a_n\otimes m + \sum_{p=0}^n (-1)^{p(n-p)} a_0\otimes ... \otimes a_p
\otimes \mu_{n-p}(a_{p+1}\otimes ... \otimes a_n\otimes m).$$

The following is straightforward.

\proclaim{(3.4.9) Proposition} (a) The differential $d$ satisfies $d^2=0$
and makes $\operatorname{ Bar}_A(M)$ into a left dg-module over 
$A$.\hfill\break
(b) The projection $\operatorname{ Bar}_A(M)\to M$ which on $A^{\otimes 
n}\otimes M$
is given by $\mu_n$, is a quasiisomorphism.\hfill\break
(c) Let $N$ be any genuine dg-module over $A$. Then an $A_\infty$-morphism
$M\to N$ is the same as a morphism of dg-modules $\operatorname{ Bar}_A(M)
\to N$.\endproclaim

\vskip .2cm

\subhead{ (3.5) A model for $\operatorname{ RAct}$ 
classifying $A_\infty$-actions}
\endsubhead Let $A$ be a $\Bbb Z_-$-graded associative dg-algebra and $V$
be an ungraded finite-dimensional vector space. We set
$$\widetilde{\operatorname{ R}} \operatorname{ Act}(A,V) = 
\operatorname{ Act}(D(A), V).\leqno (3.5.1)$$
So it is a model for $\operatorname{ RAct}(A,V)$, defined via
the particular quasifree resolution $D(A)$ of $A$. 
We postpone till n.(3.6) the discussion of other resolutions and concentrate
on this model. 
By construction, the affine dg-scheme
$\widetilde{\operatorname{ R}} \operatorname{ Act}(A,V)$ is the classifier  of 
$A_\infty$-actions.
Its coordinate ring 
$\Bbb K[\widetilde{\operatorname{ R}}\operatorname{ Act}(A,V)]$ is
the {\it free} graded commutative algebra on the
matrix elements of {\it indeterminate} linear operators $\mu_n: A^{\otimes n}
\otimes V\to V$  while the differential is chosen so as to enforce
(3.4.1). In other words, we have:

\proclaim{(3.5.2) Proposition} For any commutative dg-algebra $\Lambda$ the 
set
$$\operatorname{Hom}_{\operatorname{ dg-Alg}}
(\Bbb K[\widetilde{\operatorname{ R}}\operatorname{ Act}(A,V)], \Lambda)$$ 
is naturally
identified with the set of $\Lambda$-(multi)linear $A_\infty$-actions of
$A\otimes_\Bbb K\Lambda$ on $V\otimes_\Bbb K\Lambda$. \endproclaim

Notice that if $A$ has all its graded components finite-dimensional,
then so does $\Bbb K[\widetilde{\operatorname{ R}}\operatorname{ Act}(A,V)]$  
and therefore
$\widetilde{\operatorname{ R}}\operatorname{ Act}(A,V)$ is a dg-manifold.

\vskip .2cm

We now describe a version of the Eilenberg-Moore spectral sequences
for the functor 
$A\mapsto \Bbb K[\widetilde{\operatorname{ R}}\operatorname{ Act}(A,V)]$.

\proclaim{(3.5.3) Proposition} For any $\Bbb Z_-$-graded associative
dg-algebra $A$ we have natural convergent spectral sequences
$$E_1 = H^\bullet \,\Bbb K[\widetilde{\operatorname{ R}}
\operatorname{ Act}(A^\bullet_\# ,V)] 
\quad \Longrightarrow \quad H^\bullet \, 
\Bbb K[\widetilde{\operatorname{ R}}\operatorname{ Act}(A,V)];
\leqno (\operatorname{ a})$$
$$E_2 = H^\bullet \, \Bbb K[\widetilde{\operatorname{ R}}
\operatorname{ Act}(H^\bullet(A), V)]
\quad \Longrightarrow \quad H^\bullet \, 
\Bbb K[\widetilde{\operatorname{ R}}\operatorname{ Act}(A,V)].
\leqno (\operatorname{ b})$$\endproclaim

\demo{Proof} Let us construct the sequence (b), the first
one being similar. As a vector space, 
$$\Bbb K[\widetilde{\operatorname{ R}}\operatorname{ Act}(A,V)] \quad =\quad 
S\biggl(\bigoplus_{n=1}^\infty
\operatorname{ Hom}(V, A^{\otimes n}\otimes V)\biggr)$$
 and its grading comes from a natural bigrading of which
the first component is induced by the grading in $A$ while the other
one is the grading in the symmetric algebra induced by the grading
of the generators 
$\operatorname{ deg}\, \operatorname{ Hom}(V, A^{\otimes n}\otimes V) = 1-n$.
Similarly, the differential $d$ is a sum $d'+d''$ where $d'$, of bidegree
$(1,0)$, is induced by the differential in $A$ and $d''$, of bidegree
$(0,1)$, is induced by the algebra structure in $A$ (i.e., $d''$ 
is the differential in $\Bbb K[\widetilde{\operatorname{ R}}
\operatorname{ Act}(A^\bullet_\# ,V)]$).
Thus we have a double complex. Now, since taking cohomology
commutes with tensor products over $\Bbb K$, we find that
$$H^\bullet_{d'}\, \Bbb K[\widetilde{\operatorname{ R}}
\operatorname{ Act}(A,V)] \simeq 
\Bbb K[\widetilde{\operatorname{ R}}\operatorname{ Act}(H^\bullet(A), V)]$$
as complexes, if we take the differential on the right to be induced by $d''$. 
So our statement follows from the standard spectral sequence of
a double complex, which converges as the double complex is
$\Bbb Z_-\times \Bbb Z_-$-graded. \enddemo

\proclaim{(3.5.4) Proposition} Suppose $A$ is concentrated in degree 0.
Then: \hfill\break
(a) We have $\pi_0 \, \widetilde{\operatorname{ R}}
\operatorname{ Act}(A,V) = \operatorname{ Act}(A,V)$.
\hfill\break
(b) For any $\mu\in \operatorname{ Act}(A,V)$ the spaces $H^i T^\bullet_\mu
\widetilde{\operatorname{ R}}\operatorname{ Act}(A,V)$ 
are given by (3.3.4).\endproclaim

\demo{Proof} (a) The underlying ordinary scheme of
$\widetilde{\operatorname{ R}}\operatorname{ Act}(A,V)$ is  the affine space
$$|\operatorname{ Hom}_\Bbb K(A\otimes V, V)|.$$ The ideal of the subscheme 
$\pi_0$
is the image, under $d$, of the $(-1)$st graded component of the coordinate
ring. The space of generators of the coordinate ring in degree $(-1)$
is $\operatorname{ Hom}_\Bbb K(V, A\otimes A\otimes V)$ and the ideal
in question is exactly given by the associativity conditions (3.3.1).
\vskip .1cm

(b) A direct inspection shows that we have an identification of
complexes
$$T^\bullet_\mu \widetilde{\operatorname{ R}}\operatorname{ Act}(A,V) = 
\operatorname{ Hom}_A(\operatorname{ Bar}^{\leq -1}_A
(V), V)[1],\leqno (3.5.5)$$
so our statement follows from the fact that $\operatorname{ Bar}_A(M)$, being
a free resolution, can be used to calculate the Ext's.\enddemo

\vskip .3cm

\subhead{ (3.6) M-homotopies}\endsubhead  To prove Proposition 3.3.6,
we  need a particular nonlinear
generalization of the principle (well known in the usual homological
algebra) that any two free resolutions of a module are homotopy equivalent. 
In order for such a statement to be useful, it needs to employ a 
concept of homotopy which is preserved under functorial constructions
on algebras. The usual notion of chain homotopy of morphisms of complexes
is preserved only under additive functors and so is not good for
our purposes. A better concept of homotopy in the nonlinear context,
which we now describe, goes back to Quillen \cite{Q, Ch. 1, Def. 4}
 cf. also \cite{BG, \S 6},
\cite {Le, Ch. II, \S 1}.

Let $A,B$ be associative dg-algebras over 
$\Bbb K$ and $(f_t: A\to B)_{t\in [0,1]}$ 
be a smooth family of dg-homomorphisms
parametrized by the unit interval in $\Bbb R$. Then, for each $t$, the 
derivative
$f'_t = {d\over dt} f_t$, satisfies
$$f'_t(ab) = f_t(a) f'_t(b) + f'_t(a) f_t(b)$$
i.e., it is a degree 0 derivation $A\to B$ with respect
to the $A$-bimodule structure on $B$ given by $f_t$. 
Also, $f'_t$ commutes with the differentials in $A$ and $B$, i.e., 
$[d,f'_t]=0$.

\proclaim{(3.6.1) Definition} An M-homotopy (M for multiplicative)
is a pair $(f_t, s_t)_{t\in [0,1]}$ where $(f_t)$ is as above and
$s_t: A\to B[-1]$ is a smooth family of degree $(-1)$ derivations
(with respect to the bimodule structures given by the $f_t$)
such that $f'_t = [d, s_t]$.\endproclaim

\proclaim{(3.6.2) Proposition} For an M-homotopy, $f_0$ and $f_1$ induce
the same morphism $H(A)\to H(B)$.\endproclaim

\demo{Proof} Clear, as $f'_t$, being homotopic to 0 in the
usual sense of cochain complexes, induces 0 on the homology. \enddemo

\vskip .1cm

\noindent {\bf (3.6.3) Remark.} A polynomial M-homotopy is the same
as a morphism of dg-algebras
$$A\to B\otimes_\Bbb K\biggl( \Bbb K[t,\epsilon], \,\,\, \operatorname{ deg}\, 
t=0,
\,\,\, \operatorname{ deg}(\epsilon)=+1, \,\,\, dt = \epsilon \biggr).$$
The dg-algebra $\Bbb K[t, \epsilon]$ on the right is $\Bbb Z_+$-graded,
so it is formally outside the framework of this paper. Nevertheless,
it is quasiisomorphic
to $\Bbb K$, so from a wider derived-categorical point of view an $M$-homotopy
should be thought of as representing  one morphism $A\to B$.

\vskip .2cm

The following construction was presented by M. Kontsevich in his course
on deformation theory (Berkeley 1994).

\proclaim{(3.6.4) Proposition} Let $B,C$ be $\Bbb Z_-$-graded dg-algebras
such that $B$ is quasifree and $C$ is acyclic in degrees $<0$.
Let $f_0, f_1: 
B\to C$ be two morphisms of dg-algebras inducing the same morphism
$H^0(B)\to H^0(C)$. Then there exists a polynomial $M$-homotopy
between $f_0$ and $f_1$. \endproclaim

\demo{Proof} As  $B$ is quasifree,
let us write $B_{\# }=F(E^\bullet)$,
 for some $\Bbb Z_-$-graded vector space $E^\bullet$ of generators.
A morphism $g: B\to C$ is uniquely defined
by its restriction on the generators which
furnishes a family of linear maps $g^{(i)}: E^{-i}
\to C^{-i}$. Conversely, any choice of such maps
which is compatible with the differentials, defines
a homomorhism. Similarly, a derivation
$\sigma: B\to C$ (with respect to the bimodule structure on
$C$ given by $g$), is uniquely described by its restriction
on generators which gives
linear maps $\sigma^{(i)}: E^{-i}\to C^{-i-1}$.

We now construct  a family of homomorphisms $(f_t): B\to C$,
interpolating between $f_0, f_1$ 
inductively, by constructing successively the $f_t^{(i)}$, $i=0,1,...$.
 To start, we define $f_t^{(0)}$ by linear interpolation:
$f^{(0)}_t(e) = (1-t)f_0(e) + t f_1(e)$, $e\in E^0$.
 On this stage the compatibility with
the differential does not yet arise. By construction, the images of 
$f_t^{(0)}(e)$
in $H^0(C)$ are independent on $t$ and therefore
$(d/dt)f_t^{(0)}(e)$ takes values in $\operatorname{ Im}\{d: C^{-1}\to C^0\}$.
So we can find a polynomial family of maps $s_t^{(0)}: E^0_1\to
C^{-1}$ such that $(d/dt)f_t^{(0)}(e) = d s_t^{(0)}(e)$. 
To continue, we need to define $f_t^{(1)}: E^{-1}\to C^{-1}$
in such a way that 
$$df_t^{(1)}(e) = f_t^{(0)}(de)\leqno (3.6.5)$$
But $f_t^{(0)}(de)$, being a linear interpolation between
$f_0(de) = df_0(e)$ and $f_1(de)=df_1(e)$, lies, for any $t$,
in the image of $d$. Therefore we can choose a polynomial family
$(f_t^{(1)})$ interpolating between $f_0^{(1)}$ and $f_1^{(1)}$
and satisfying (3.6.5). Next, for $e\in E^{-1}$, we have
$$d\biggl( {d\over dt} f_t^{(1)}(e) - s^{(0}_t (de)\biggr) = 
d\biggl( {d\over dt} f_t^{(0)}(de)\biggr)-d(s^{(0)}_t(de)) =
d(s^{(0)}_t(de))-d(s^{(0)}_t(de))=0$$
and because $B_2$ is acyclic in degrees $\leq -1$, we can find
a polynomial family of linear maps $(s_t^{(1)}: E_1^{-1}\to B_2^{-2})$
such that
$$d(s_t^{(1)}(e)) = {d\over dt} f_t^{(1)}(e) - s_t^{(0)}(de)$$
which is the first in the series of conditions defining an M-homotopy.
We then continue in this way, defining successively the $f_t^{(i)}$
and $s_t^{(i)}$ on $E^{-i}$ and extending them to homomorphisms
(resp. derivations) on the subalgebra generated by $E^{-i}, ..., E^0$.
This furnishes a required M-homotopy.\enddemo

\vskip .2cm

Let now $A=F(E^\bullet)$ be the free associative algebra on the
$\Bbb Z_-$-graded vector space $E^\bullet$ (no differential). Then the
quasiisomorphism $\alpha: D(A)\to A$ described in (3.4.5)(b),
has a natural
right inverse 
$\beta: A\to D(A)$, so that $\alpha\beta = \operatorname{ Id}_A$. More 
precisely, $\beta$ is defined on the space of generators $E\subset A$ to
identify it with the natural copy of $E$ inside $A\subset F(A)\subset D(A)$
and then extended to the entire $A$ because $A$ is free. 
Thus $\beta$ is also a quasiisomorphism.

\proclaim{(3.6.6) Proposition} The composition $\beta\alpha: D(A)\to D(A)$
is M-homotopic to the identity of $D(A)$. 
\endproclaim

\demo{ Proof} If $E^\bullet$ is in degree 0, then so is $A$
and thus we can apply Proposition 3.6.4 to $B=C=D(A)$. If $E^\bullet$
is not concentrated in degree 0, then we notice that $D(A)$
in fact comes from a $\Bbb Z_-\times \Bbb Z_-$-graded dg-algebra
and that we can mimic all the steps in the proof of (3.6.4), using
the induction in the second component of the bidegree.\enddemo

\vskip .3cm

\subhead {(3.7) Proof of Proposition 3.3.6}\endsubhead 
Proposition 3.6.6 implies the following.

\proclaim{(3.7.1) Corollary} If $A = F(E^\bullet)$ is free with trivial 
differential, then $\widetilde{\operatorname{ R}}\operatorname{ Act}(A,V)$ is 
quasiisomorphic to
$\operatorname{ Act}(A,V) = 
|\operatorname{ Hom}_\Bbb K(E^\bullet\otimes V, V)|$.\endproclaim

\demo{Proof} As we pointed out before, the correspondence
$A\mapsto \operatorname{ Act}(A,V)$ is contravariantly functorial
in $A$; equivalently, $\Bbb K[\operatorname{ Act}(A,V)]$ depends on $A$
in a covariant way. Thus the maps
 $\alpha,\beta$ between $A$ and $D(A)$ induce morphisms of
commutative dg-algebras $\alpha_*, \beta_*$ from
$\Bbb K[\operatorname{ Act}(A,V)]$ to $\Bbb K[\operatorname{ Act}(D(A), V)]=
\Bbb K[\widetilde{\operatorname{ R}}\operatorname{ Act}(A,V)]$ and back, with 
$\alpha_*\beta_* = \operatorname{
Id}$. Further, the polynomial M-homotopy between  
$\beta\alpha$ and Id, constructed in (3.6.6), is also inherited,
because of Remark 3.6.3,
in functorial constructions such as passing to $\Bbb K[\operatorname{ Act}(-, 
V)]$.
This proves the statement. 
\enddemo

\proclaim{(3.7.2) Proposition} If $B\buildrel p\over\to A$ is any quasifree
associative dg-resolution, then $\Bbb K[\operatorname{ Act}(B,V)]$
is naturally quasiisomorphic to $\Bbb K[\widetilde{\operatorname{ 
R}}\operatorname{ Act}(A,V)]$
and therefore it is independent, up to a quasiisomorphism,
of the choice of $B$. \endproclaim

\demo{Proof} 
By Corollary 3.6.5, $\Bbb K[\operatorname{ Act}(B_\# ,V)]$
is quasiisomorphic to $\Bbb K[\widetilde{\operatorname{ R}}\operatorname{ 
Act}(B_\# ,V)]$,
the quasiisomorphism being induced by the map $\alpha$.
Proposition 3.5.3(a) implies then that
the map 
$$\Bbb K[\operatorname{ Act}(B,V)] \to
\Bbb K[\widetilde{\operatorname{ R}}\operatorname{ Act}(B,V)]
$$
is also a quasiisomorphism as it induces an isomorphism of the first
terms of the spectral sequences described in 3.5.3(a). 
Further, the spectral sequence (3.5.3)(b) shows that the morphism
$$p_*: \Bbb K[\widetilde{\operatorname{ R}}\operatorname{ Act}(B,V)] \to
\Bbb K[\widetilde{\operatorname{ R}}\operatorname{ Act}(A,V)]$$
is a quasiisomorphism. This proves our statement.\enddemo

\vskip .1cm

Thus we have established Proposition 3.3.6.

\vskip .3cm

\subhead{(3.8) The derived linearity locus}\endsubhead Let ${ S}$ be a
$\Bbb Z_-$-graded dg-scheme and $A$ be a $\Bbb Z_-$-graded
associative dg-algebra. Let $M,N$ be two quasicoherent dg-sheaves on $S$ such 
that
$M_\# , N_\# $ are locally free over ${\Cal   O}_{S^\# }$.
We assume that the generators of $M_\# $ are in degrees $\leq 0$
and those of $N_\# $ are in degrees $\geq 0$. 
Suppose that $M,N$ are made into dg-modules over $A\otimes_\Bbb K{\Cal   O}_S$
and we have a morphism $f: M\to N$ of ${\Cal   O}_S$-dg-modules
(but not necessarily of $A\otimes_\Bbb K{\Cal   O}_S$-dg-modules).
According to the general approach of homological
algebra, we define the derived linearity locus $\operatorname{ RLin}_A(f)$ as
the derived fiber product (2.8)
$$\minCDarrowwidth{7 mm}
\CD \operatorname{ RLin}_A(f)@>>> S\\
@VVV @VV f V\\
|\operatorname{ R} \Cal Hom_{A\otimes{\Cal   O}_S}(M,N)| @>|\rho| >> 
|\Cal Hom_{{\Cal   O}_S}(M,N)|.\endCD\leqno (3.8.1)$$
Here $\operatorname{ R} \Cal Hom_{A\otimes{\Cal   O}_S}(M,N) = 
\Cal Hom_{A\otimes{\Cal   O}_S}(P,N)$, where
$P\to M$ is a resolution by a dg-module such that:

\vskip .1cm

\item{(1)} $P$ is $\Bbb Z_-$-graded and
$P_\# $ is locally (on the Zariski topology of $S$)
projective over $A_\# \otimes {\Cal   O}_{S^\# }$. 

\vskip .1cm

 If, in addition, we have a stronger condition, namely:

\vskip .1cm

\item{(2)} The  morphism $\rho$ is
a termwise surjective morphism of cochain complexes.

\vskip .1cm

\noindent then the derived fiber product coincides with
the usual fiber product, as the morphism $|\rho|$ is flat.  

\vskip .1cm

One example of a resolution satisfying (1) and (2) above is
the bar-resolution $\operatorname{ Bar}_A(M)\to M$,
see (3.4.7) (we consider $M$ as an $A_\infty$-module with $\mu_i=0, i\geq 2$).
The following is then a standard application of the Eilenberg-Moore
spectral sequences for the derived Hom and tensor product
of dg-modules.

\proclaim{(3.8.2) Proposition}  The derived linearity locus is
independent,
up to a quasiisomorphism, on the choice of  $P$
satisfying (1) and (2).\endproclaim

 We  denote by 
$$\widetilde{\operatorname{ R} }
\operatorname{ Lin}_A(f) =|\Cal Hom^\bullet_{A\otimes 
{\Cal   O}_S}(\operatorname{ Bar}_A(M),N)|
\times_{|\Cal Hom_{{\Cal   O}_S}(M,N)|}
S \leqno (3.8.3)$$ 
the particular model for RLin
obtained by using the bar-resolution.
 (we consider $M$ as an $A_\infty$-module with $\mu_i=0, i\geq 2$).
Let us note some additional properties of this model.
First, it can be applied in a more general situation. Namely,
let $A$ and $N$ be as before, but assume that $M$ is only
an $A_\infty$-module over $A$. In this case, as $\operatorname{ Bar}_A(M)$ 
makes sense,
we define  $\widetilde{\operatorname{ R}}\operatorname{ Lin}_A(f)$ by (3.8.3). 
 Notice that while we can view
an $A_\infty$-morphism $f: M\to N$ as a morphism of $D(A)$-dg-modules,
(still denoted $f$), the $A_\infty$-version $\widetilde{\operatorname{ 
R}}\operatorname{ Lin}_A(f)$
is much more economical in size than any of the models for
$\operatorname{ RLin}_{D(A)}(f)$ given by (3.8.1-2), especially than
$\widetilde{\operatorname{ R}}\operatorname{ Lin}_{D(A)}(f)$.

\proclaim{(3.8.4) Proposition} 
$\widetilde{\operatorname{ R}}\operatorname{ Lin}_A(f)$
is the complex of vector bundles on $S$
$$\operatorname{ Cone} \biggl\{ {\Cal   O}_S
\buildrel\delta f\over\longrightarrow
\Cal Hom_{{\Cal   O}_S}
(\operatorname{ Bar}_A(M), N)\biggr\}[1]$$
considered as a dg-scheme.
Here $\delta f\in \Cal Hom_{{\Cal   O}_S}(A\otimes M, N)$
takes $a\otimes m \mapsto f(a\otimes m)-af(m)$.\endproclaim

As with $\widetilde{\operatorname{ R}}\operatorname{ Act}$, the  construction 
of
$\widetilde{\operatorname{ R}}\operatorname{ Lin}$ can be interpreted via
 $A_\infty$-structures.

\proclaim{(3.8.5) Proposition} (a) The natural morphism 
$p: \widetilde{\operatorname{ R}}\operatorname{ 
Lin}_A(f) \to S$ is smooth and the induced morphism $p^*f: p^*M\to p^*N$
is an $A_\infty$-morphism of dg-modules over
$A\otimes {\Cal   O}_{\widetilde{\operatorname{ R}}
\operatorname{ Lin}_A(f)}$.\hfill\break
(b)  For any commutative dg-algebra $\Lambda$ the set
$\operatorname{ Hom}_{\operatorname{ dg-Sch}}(\operatorname{ Spec}(\Lambda), 
\widetilde{\operatorname{ R}}\operatorname{ Lin}_A(f))$
is identified with the set of data $(g, h_1, h_2, ...)$
where $g: \operatorname{ Spec}(\Lambda)\to S$ is a morphism of dg-schemes and
$h_n: A^{\otimes n}\otimes_\Bbb K g^*M\to g^*N$ are such
that $(g^*f, h_1, h_2, ...)$ is an $A_\infty$-morphism $g^*M\to g^*N$.
\endproclaim

Informally,  $\widetilde{\operatorname{ R}}\operatorname{ Lin}_A(f)$ is 
obtained by
adding to ${\Cal   O}_S$ new free generators which are matrix
elements of interdeterminate higher homotopies $h_i: A^{\otimes i}\otimes M\to 
N$, $i\geq 1$ and arranging the differential there so as to satisfy (3.4.2).

\vskip .3cm

\subhead{ (3.9) The derived $A$-Grassmannian}\endsubhead
We place ourselves in the situation of the beginning of this section,
So $A$ is a finite-dimensional $\Bbb K$-algebra and $M$ a finite-dimensional
$A$-module. By applying the construction of the derived space of
actions to any fiber of the tautological bundle $\widetilde V$ on $G(k,M)$,
we get a dg-scheme 
$\operatorname{ RAct}(A,\widetilde V)\buildrel q\over\rightarrow G(k,M)$.
If we take a quasifree resolution $B\to A$ with finitely many
generators in each degree, then $\operatorname{ RAct}(A,\widetilde V)$ will be 
a dg-manifold.
For example, the model $\widetilde{\operatorname{ R}}\operatorname{ Act}(A, 
\widetilde V)$
obtained via the bar-resolution $D(A)$, satisfies this property. 
Thus $q^*\widetilde V$ is  a dg-module
over $B$.

\proclaim{(3.9.1) Definition} The derived $A$-Grassmannian $RG_A(k,M)$
is defined as the derived linearity locus $\operatorname{ RLin}_{B}(f)$,
where $f: q^*\widetilde V\to \underline{M}$ is the tautological
morphism from (3.2.3).

\endproclaim

A smaller model can be obtained by taking $B=D(A)$, viewing
a dg-module over $D(A)$ as an $A_\infty$-module over
$A$ and apply the construction of the derived linearity locus for
$A_\infty$-modules described in (3.8). This model is a dg-manifold.

\proclaim{(3.9.2) Theorem} (a) We have $\pi_0\, RG_A(k,M) = G_A(k,M)$.
\hfill\break
(b) For any $A$-submodule $V\subset M$ with $\dim_\Bbb K(M)=k$, we have
$$H^i T^\bullet_{[V]} RG_A(k,M) = \operatorname{ Ext}^i_A(V, M/V).$$
\endproclaim

\demo{Proof} (a) follows from similar properties of
RAct, RLin (Propositions 3.5.4 and 3.8.4).
To see (b), notice that we have an identification
in the derived category:
$$\operatorname{ RHom}_A(V, M/V) \sim \operatorname{ Cone}\biggl\{ 
\operatorname{ RHom}_A(V,V)\to \operatorname{ RHom}_A(V,M)\biggr\}[1].$$
To be specific, we will consider the model for $RG_A$ obtained
by using  $\widetilde{\operatorname{ R}}\operatorname{ Act}$ and
the $A_\infty$-version of $\widetilde{\operatorname{ R}}\operatorname{ Lin}$. 
Then, denoting $\mu:A\otimes V\to V$ the induced $A$-action on
the submodule $V$, we have, by Proposition 3.8.4, an identification of 
complexes
$$T^\bullet_{[V]}RG_A(k,M) = 
\operatorname{ Cone}\biggl\{ T^\bullet_{([V], \mu)}
\widetilde{\operatorname{ R}}\operatorname{ Act}(A, \widetilde V)
\to T^\bullet_{([V], \mu)}
|\Cal Hom_{{\Cal   O}_{\widetilde{\operatorname{ R}}
\operatorname{ Act}(A, \widetilde V)}}
(\operatorname{ Bar}_A(q^*\widetilde V), \underline {M})|\biggr\}.$$
The dg-scheme $\widetilde{\operatorname{ R}}\operatorname{ Act}(A, \widetilde 
V)$
is a fibration over $G(k,M)$, and  
$|\Cal Hom_{{\Cal   O}_{\widetilde{\operatorname{ R}}
\operatorname{ Act}(A, \widetilde V)}}
(\operatorname{ Bar}_A(q^*\widetilde V), \underline {M})|$ is a fibration over
 $\widetilde{\operatorname{ R}}\operatorname{ Act}(A, \widetilde V)$ 
so it also a fibration over
$G(k,M)$. The tangent bundle of each of these fibrations fits into
short exact sequence involving the relative tangent bundle and 
the pullback of $TG(k,M)$. Let us write the corresponding exact sequences
for fibers of the tangent bundles. Using Propositions 3.5.4 and 3.8.4,
we can write them as follows:
$$0\to \operatorname{ Hom}_A(\operatorname{ Bar}^{\leq -1}_A(V),V) \to
T^\bullet_{([V],\mu)}\widetilde{\operatorname{ R}}
\operatorname{ Act}(A, \widetilde V)\to
T_{[V]}G(k,M)\to 0,$$
$$0\to \operatorname{ Hom}_{A}(\operatorname{ Bar}_A(V), M)\to 
T^\bullet_{([V], \mu)}
|\Cal Hom_{{\Cal   O}_{\widetilde{\operatorname{ 
R}}\operatorname{ Act}(A, \widetilde V)}}
(\operatorname{ Bar}_A(q^*\widetilde V), \underline {M})| \to 
T_{[V]}G(k,M)\to 0.$$
This means that in the cone the two copies of $T_{[V]}G(k,M)\to 0$
will cancel out, up to quasiisomorphism, and we conclude that
$$T_{[V]}RG_A(k,M) = \operatorname{ Cone}\biggl\{\operatorname{ 
Hom}_A(\operatorname{ Bar}_A^{\leq -1}(V),
V)\to \operatorname{ Hom}_A(\operatorname{ Bar}_A(V),M)\biggr\} [1],$$
whence the statement. \enddemo

\vskip .2cm

\noindent {\bf (3.9.3) Remarks.} (a) Instead of working with the
module structures on the fibers of the universal subbundle $\tilde V$
on $G(k, M)$, we could equally well work with module structures
on the fibers of the universal quotient bundle $M/\tilde V$ and modify the
approach of (3.2) accordingly. 

(b) Taking $M=A$, set $\Cal J(k, A)= G_A(k, A)$ (the scheme of
ideals in $A$ of dimension $k$). When $A$ is commutative,
there is a derived analog $R\Cal J(k, A)$ of $\Cal J(k, A)$
different from $RG_A(k, A)$ and  whose
 construction will be described in detail in \cite{CK}. 
The approach is based on realizing $\Cal J(k, A)$ via
two constructions similar but not identical to those
described in (3.2). The first one is the space  $\Cal C(W)
\subset \text{Hom}(S^2W, W)$ of all
commutative algebra structures on a finite-dimensional vector space $W$.
Applying this to fibers of the bundle $A/\tilde V$ on $G(k, A)$,
we get a fibration $q: \Cal C(A/\tilde V)\to G(k, V)$ and a 
vector bundle morphism
$g: \underline{A}\to q^*(A/\tilde V)$ on $\Cal C(A/\tilde V)$.
Fibers of both these bundles are commutative algebras and 
  $\Cal J(k, A)$
is the {\it homomorphicity locus} of $g$, i.e., the subscheme of points of
the base such that the corresponding morphism of the fibers is an
algebra homomorphism. The dg-manifold $R\Cal J(k, A)$  is obtained
by taking the derived versions of these steps. It
will be used in constructing the derived Hilbert
scheme mentioned in (0.4).

\vfill\eject

\heading{4. Derived Quot schemes}\endheading

In this section we will apply the construction of $RG_A(k,M)$
of \S 3 to the case of interest in geometry, when $A = \bigoplus_{i\geq 0}
 H^0(X, \Cal  O_X(i))$ for a projective scheme $X$ 
and $M=\bigoplus_{i=p}^q H^0(X, \Cal  F(i))$ for a coherent
sheaf $\Cal  F$ on $X$.  In this situation all objects acquire
extra grading and to avoid confusion, we sharpen our terminology.

\subhead{(4.1) Conventions on grading}\endsubhead
We will consider  bigraded vector spaces $V = V^\bullet_\bullet = 
\bigoplus_{p,q} V^p_q$.
The lower grading will be called {\it projective} and the upper
one, {\it cohomological}. By a bigraded complex we mean a bigraded
vector space with a differential $d$ having degree 1 in the upper
grading and 0 in the lower one. Tensor products $V^\bullet_\bullet\otimes 
W^\bullet_\bullet$ of bigraded complexes
are defined in the usual way and the symmetry map 
$V^\bullet_\bullet\otimes W^\bullet_\bullet \to 
W^\bullet_\bullet\otimes V^\bullet_\bullet$ is defined by the Koszul
sign rule involving only the upper grading. 
The concepts of a bigraded (commutative) dg-algebra, bigraded $A_\infty$-
algebra etc. will be understood accordingly, with only the upper
grading contributing to the sign factors.

Given a (lower) graded associative algebra $A=\bigoplus A_i$
and its left graded modules $M = \bigoplus M_i$, $N=\bigoplus N_i$
we define the $\operatorname {Ext}^{i,0}_A(M,N)$ to be the
derived functors of $\operatorname{Hom}^0_A(M,N)$, i.e., of
the Hom functor in the category of graded modules.

\subhead{(4.2) Derived $A$-Grassmannian in the graded case}\endsubhead
Let $A=\bigoplus_{i\geq 0}A_i$ be a graded associative
algebra with $A_0=\Bbb K$ and $\dim(A_i)<\infty$ for all $i$.
Let $M = \bigoplus M_i$ be a finite-dimensional graded $A$-module;
thus there are $p\leq q$ such that $M_i=0$ unless $i\in [p,q]$. 
If $k=(k_i)$ is a sequence of nonnegative integers, we have introduced
in (1.3) the graded $A$-Grassmannian $G_A(k,M)$. The construction
of the derived $A$-Grassmannian $RG_A(k,M)$ in Section 3, can be
repeated without changes in the graded case, if we consider
everywhere only morphisms of degree 0 with respect to the
lower grading and replace the functor Hom with $\text{Hom}^0$.

\vskip .2cm

\proclaim{ (4.2.1) Theorem} 
 Up to quasiisomorphism, the graded version of the derived $A$-Grassmannian is 
a dg-manifold
satisfying the conditions:
$$\pi_0 RG_A(k, M) = G_A(k, M), \quad H^i T_{[V]}  RG_A(k, M) = 
\operatorname{Ext}^{i,0}_A(V, M/V).$$
\endproclaim

\demo{Proof} The only issue that needs to be addressed is the existence of
a model for $RG_A(k,M)$ which is of finite type, as $A$ is now 
infinite-dimensional over $\Bbb K$. More precisely, we need to show that
the particular model $\widetilde{\text {R}}\text{Act}(A, V)$ and
$\widetilde {\text{R}}\text{Lin}_A(f)$ are of finite type
in the bigraded context as well 
(then they will be dg-manifolds by construction).

To see this, recall that  $\widetilde{\text {R}}\text{Act}(A, V)$
is the affine dg-scheme whose coordinate algebra 
$\Bbb K [\widetilde{\text {R}}\text{Act}(A, V)]$ is the free (upper)
graded commutative algebra on the matrix elements of
indeterminate linear maps $A^{\otimes n}\otimes V\to V$ 
{\it of degree 0 with respect to the lower grading}. Since $V$
(a graded subspace of $M$)
is concentrated in only finitely many degrees (from $p$ to $q$),
 and since we
can disregard $A_0=\Bbb K$, there are only finitely many possibilities
for nonzero maps $A_{i_1}\otimes ... \otimes A_{i_n}\otimes V_j
\to V_{i_1+...+i_n+j}$, $i_\nu > 0$, $p\leq j\leq q$. Each such
possibility gives a finite-dimensional space of maps. This
implies that each (upper) graded component of 
$\Bbb K [\widetilde{\text {R}}\text{Act}(A, V)]$ is finite-dimensional, so
$\widetilde{\text {R}}\text{Act}(A, V)$ is a dg-manifold.  
Proposition 3.5.4 now holds in the graded context, with Hom and
Ext replaced everywhere by $\text{Hom}^0$ and $\text{Ext}^0$. 

Further, if we use the same convention in (3.8), we get that 
the graded version $\widetilde{\text{R}}
\text{Lin}$ is also a dg-manifold. The theorem is proved.

\enddemo

 \subhead{(4.3) The derived Quot scheme}\endsubhead
Let now $X\subset\Bbb P^n$ be a {\it smooth} projective variety, $\Cal  F$ a
coherent sheaf on $X$ and $h\in\Bbb Q[t]$ a polynomial. Let $A$ be the graded
coordinate algebra of $X$ and $M$ the graded $A$-module corresponding to 
$\Cal  F$, see (1.2). 

\proclaim{(4.3.1) Definition} The derived Quot scheme is defined as
$$RSub_h(\Cal  F):=RG_A(h,M_{[p,q]})\ \ \text{for}\
0\ll p\ll q.$$
Here $RG_A(h,M_{[p,q]})$ is the graded version of the derived
Grassmannian constructed in (4.2).
\endproclaim

The well-definedness of $RSub_h(\Cal  F)$ up to isomorphism in the derived
category $\Cal  D Sch$
is part (a) of the following theorem which is the main result of this
paper. 

\proclaim{(4.3.2) Theorem} (a) For $0\ll p\ll p'\ll q'\ll q$ the natural
projection
$$RG_A(h, M_{[p,q]})\to RG_A(h, M_{[p', q']})$$
is a quasiisomorphism of dg-manifolds.

(b) $\pi_0 RSub_h(\Cal  F) = Sub_h(\Cal  F)$.

(c) If $\Cal  K\subset\Cal  F$ has Hilbert polynomial $h$, then
$$H^iT^\bullet_{[\Cal  K]} RSub_h(\Cal  F) \simeq \operatorname
{Ext}^i_{\Cal  O_X}(\Cal  K, \Cal  F/\Cal  K), \quad i\geq 0.$$
\endproclaim

\demo{Proof}
Part (b) follows from Theorem 1.4.1 and (4.2.1). Part (a) would follow
from (b)(c) in virtue of the ``Whitehead theorem" (2.5.9). 
More precisely, we need to apply (c) to the dg-scheme
$RSub_h(\Cal F)\otimes \Bbb F$ for any field extension
$\Bbb F$ of $\Bbb K$. This scheme is just the $RSub$ scheme
corresponding to the sheaf $\Cal F\otimes \Bbb F$ on the
scheme $X\otimes\Bbb F$.
So we concentrate
on (c) and start with the following.
\enddemo

\proclaim{(4.3.3) Proposition} (a) If $\Cal  F, \Cal  G$  are coherent
sheaves on $X$ with corresponding graded $A$-modules
$M=\text{Mod}(\Cal  F)$ and $N=\text{Mod}(\Cal  G)$ respectively, then
$$\operatorname{Ext}_{\Cal  O_X}^i(\Cal  F, \Cal  G)=
\lim_{\buildrel p\over\rightarrow}\operatorname{Ext}_A^{i,0}(M_{\geq p},
N_{\geq p}),$$
and the limit is achieved.

(b) There exists an integer $p$ such that
$$\operatorname{Ext}_{\Cal  O_X}^i(\Cal  K, \Cal  F/\Cal  K)=
\operatorname{Ext}_A^{i,0}(W_{\geq p},
M_{\geq p}/W_{\geq p}),\quad W=\operatorname{Mod}(\Cal  K),$$
for all subsheaves $\Cal  K$ of $\Cal  F$ representing $\Bbb K$-points of
$Sub_h(\Cal  F)$, where $W$ is the graded $A$-module corresponding
to $\Cal  K$ and $M$ is the graded $A$-module corresponding
to $\Cal  F$.
\endproclaim

\demo{Proof} Part (a) follows from Serre's theorem (1.2.2); part (b)
follows from (a), from semicontinuity of the rank of a matrix and
from the fact that $Sub_h(\Cal  F)$ is a scheme of finite type.
\enddemo

\vskip .1cm

We now continue the proof of (4.3.2)(c). 
 Since $X$ is smooth, $\operatorname{Ext}_A^{i,0}\neq 0$ for only
finitely many $i$'s. In view of Proposition 4.3.3  we are reduced to

\proclaim{(4.3.4) Proposition} Let $M$, $N$, be any finitely generated graded
$A$-modules. Then for any fixed $i$ there exists an integer $q_0$ such
that
$$
\operatorname{Ext}_A^{i,0}(M,N)=
\operatorname{Ext}_A^{i,0}(M_{\leq q},N_{\leq q})
$$ 
for all $q\geq q_0$. Moreover, if $M_s$ and $N_s$ vary in a family parametrized
by a projective scheme $S$,
 then $q_0$ can be chosen independent on $s\in S$.
\endproclaim
\demo{Proof} Assume first that $M$ is free, i.e. $M=A\otimes_{\Bbb 
C}E_\bullet$,
with $E_\bullet$ a finite dimensional graded $\Bbb K$-vector space. 
If $i=0$, we have obviously
$$\operatorname{Hom}_A^0(M_{\leq q},N_{\leq q})=\operatorname{Hom}_A^0(M,N)
=\operatorname{Hom}_{\Bbb K}^0(E_\bullet,N) \leqno (4.3.5)$$
whenever $q$ exceeds the maximum of the degrees of the 
nonzero graded components
of $E_\bullet$.

Next we claim that (for $M$ free)
$$ \operatorname{Ext}_A^{i,0}(M_{\leq q},N_{\leq q})=0,\ \text{for\ all}\
i>0,\ \text{all}\ q\geq 0\ \text{and\ any}\ N.\leqno (4.3.6)
$$ 
Indeed, the long exact sequence
$$0\longrightarrow M_{\geq q+1}\longrightarrow M\longrightarrow M_{\leq q}
\longrightarrow 0$$
induces (for $i>0$) an exact sequence
$$\operatorname{Ext}_A^{i-1,0}(M_{\geq q+1},N_{\leq q})\longrightarrow
\operatorname{Ext}_A^{i,0}(M_{\leq q},N_{\leq q})\longrightarrow
\operatorname{Ext}_A^{i,0}(M,N_{\leq q}).$$
But the last term in the above sequence vanishes, since $M$ is free. Hence
(4.3.6) follows from:
 
\proclaim{(4.3.7) Lemma} For any graded $A$-modules $M$ and $N$ and all
$i\geq 0$, we have $\operatorname{Ext}_A^{i,0}(M_{\geq q+1},N_{\leq q})=0$.
\endproclaim

\demo{Proof of Lemma (4.3.7)}
For $i=0$ the statement is obvious, as $M_{\geq q+1}$ and $N_{\leq q}$ have
nontrivial graded components only in disjoint ranges of degrees. For $i>0$,
the 
groups $\operatorname{Ext}_A^{i,0}(M_{\geq q+1},N_{\leq q})$ are calculated 
as the cohomology of the complex $\operatorname{Hom}_A^0(P^{\bullet},
N_{\leq q})$,
where $P^{\bullet}$ is a free homogeneous resolution of $M_{\geq q+1}$. 
This resolution can be chosen such that each $P^j$ is concentrated in
degrees (with respect to the lower grading) at least $q+1$, so 
$\operatorname{Hom}_A^0(P^j,N_{\leq q})=0$
for all $j$ and the lemma follows.
\enddemo

Proposition  4.3.4 is therefore true for $M$ free. If $M$ is now arbitrary,
let
 $$F^{\bullet}=\{\cdots\longrightarrow F^{-1}\longrightarrow
F^0\}\longrightarrow M$$
 be a free resolution
with $F^{-j}=A\otimes_{\Bbb K} E^{-j}_\bullet$ and each $E^{-j}_\bullet$ a 
finite
dimensional graded (by lower grading) vector space. The truncation 
$F^{\bullet}_{\leq q}$
is then a resolution of $M_{\leq q}$, but it is not free anymore.
However, by (4.3.6) and the ``abstract De Rham theorem''
$\operatorname{Ext}_A^{i,0}(M_{\leq q},N_{\leq q})$ 
is still calculated by the cohomology of the complex 
$\operatorname{Hom}_A^0(F^{\bullet}_{\leq q},N_{\leq q})$. 
By (4.3.5), for some
fixed $i$, the first $i+1$ terms of this complex will be the same as the
first $i+1$ terms of $\operatorname{Hom}_A^0(F^{\bullet},N)$ whenever
$q\geq$ maximum of the degrees of the nonzero graded components
of all $E^{-j}$, $0\leq j\leq i+1$. Hence
$$
\operatorname{Ext}_A^{j,0}(M,N)=
\operatorname{Ext}_A^{j,0}(M_{\leq q},N_{\leq q})
$$
for all $j\leq i$. Finally, notice that the above proof also shows
the existence of a lower bound for $q_0$ in a  family of modules
parametrized by any projective scheme $S$. Such a family of modules
is just a graded $A\otimes_{\Bbb K}\Cal  O_S$-module $\Cal  M$
with graded components being locally free of finite rank over $\Cal  O_S$.
Because $S$ is projective, we can find a resolution of $\Cal  M$
by $A\otimes \Cal  O_S$-modules of the form 
$\Cal  F^j = A\otimes \Cal  E^j_\bullet$ where $\Cal  E^j_\bullet$ is a graded
vector bundle on $S$ such that $\bigoplus_p \Cal  E^i_p$ has finite
rank, and then the above arguments apply word by word. 
\enddemo 

This concludes the proof of Proposition 4.3.4 and  of Theorem 4.3.2.

\vskip .2cm

\noindent {\bf (4.3.8) Remark.}  If an algebraic group $G$ acts on
$\Cal F$ by automorphisms, then we have an induced action on
$Sub_h(\Cal F)$. The above construction of $RSub_h(\Cal F)$
via the derived $A$-Grassmannian and the model for the latter
via the bar-resolution (4.2.1) immediately imply that $G$ acts
on $RSub_h(\Cal F)$ by automorphisms of dg-manifolds. A case
particularly important for constructing the derived moduli stack
of vector bundles on $X$ is $\Cal F = \Cal O_X(-N)^{\oplus r}$,
$N, r\gg 0$ and $G=GL_r$. The $G$-action on an
appropriate open part of $RSub$ gives rise to a groupoid
in the category of dg-manifolds, and such groupoids
provide, as  $N, r\to\infty$, more and more representative charts for
the  moduli (dg-)stack. The exact way of gluing such charts
(by quasiisomorphisms) into a global dg-stack requires a
 separate treatment.

\subhead{(4.4) Independence of $RSub(\Cal F)$ on the projective embedding} 
\endsubhead 
 Clearly, the concept of the Hilbert polynomial of a coherent sheaf on
 $X$ depends on the choice of a very ample line bundle $L$
(the pullback of $\Cal O(1)$ under the projective embedding). Accordingly,
the  scheme $Sub_h(\Cal F)$, see (1.1.1) depends on the choice of $L$.
To emphasize this dependence, let us denote it $Sub_h^L(\Cal F)$. 
 It is well known, however, that the union
$$Sub(\Cal F) = \coprod_{h\in \Bbb Q[t]} Sub^L_h(\Cal F)$$
  depends on $X$ and $\Cal F$ but not
on $L$, as it parametrizes subsheaves in $\Cal F$ with flat quotients. 
The analog of this classical statement is the following fact.

\proclaim{(4.4.1) Theorem}
Let $L_1,  L_2$ be two very ample line bundles on $X$. Then we have
an isomorphism
$$\coprod_{h\in \Bbb Q[t]}RSub_h^{L_1}(\Cal F) \quad \sim \quad
\coprod_{k\in \Bbb Q[t]}RSub_k^{L_2}(\Cal F)$$
in the derived category of  dg-schemes (of infinite type).
 \endproclaim

Thus we have a well defined, up to quasiisomorphism, dg-scheme which we
can denote
$RSub(\Cal F)$.

\demo{Proof}
We begin with the following simple general fact. Let $X$ be a dg-manifold
and let $Z$ be a connected component of $\pi_0(X)$. Choose an open subscheme
$Y^0\subset X^0$ such that $Y^0\cap\pi_0(X)=Z$ and set $\Cal O_Y^{\bullet}=
\Cal O_X^{\bullet}\mid_{Y^0}$. Then $Y$ is a dg-manifold. Moreover, if we
have $Y_1$ and $Y_2$ as above, then they are quasiisomorphic to $Y_1\cap Y_2$.
The discussion above implies the following.
\enddemo

\proclaim{(4.4.2) Proposition} 
$X$ is quasiisomorphic to a dg-manifold which is a disjoint
union of open submanifolds, each of them containing exactly one 
connected component of $\pi_0(X)$.
\endproclaim 

To continue the proof of (4.4.1),   put
$$A:=\bigoplus_{m\geq 0}H^0(X,L_1^{\otimes m}),\;\; 
B:=\bigoplus_{n\geq 0}H^0(X,L_2^{\otimes n}),\;\;
C:=\bigoplus_{m,n\geq 0}H^0(X,L_1^{\otimes m}\otimes L_2^{\otimes n}),$$
and
$$M:=\bigoplus_{m\in\Bbb Z}H^0(X,\Cal F\otimes L_1^{\otimes m}),\;\;
N:=\bigoplus_{n\in\Bbb Z}H^0(X,\Cal F\otimes L_2^{\otimes n}),$$
$$P:=\bigoplus_{m,n\in\Bbb Z}H^0(X,\Cal F\otimes L_1^{\otimes m}
\otimes L_2^{\otimes n}).$$
Thus $C$ is a bigraded algebra and $P$ is a bigraded $C$-module.
Finite-dimensional truncations of $P$ will be denoted by
  $$P_{[(p_1,p_2),(q_1,q_2)]} = \bigoplus_{p_1\leq m\leq q_1, p_2\leq n\leq 
q_2}
H^0(X,\Cal F\otimes L_1^{\otimes m}
\otimes L_2^{\otimes n}).$$

The choice of $L_1, L_2$ allows one to associate to any coherent
sheaf  $\Cal G$ on  $X$ its Hilbert polynomial $\frak H^{\Cal G}(t,s)$
depending on two variables:
$$\frak H^{\Cal G}(t,s) = \operatorname{dim}
H^0(X,\Cal G\otimes L_1^{\otimes m}
\otimes L_2^{\otimes n}) \;\;\text{ for}\;\; m,n\gg 0.$$

  Let $\frak H$ be a polynomial in $\Bbb Q[s,t]$
and let $Sub_{\frak H}(\Cal F)$ be the part of the $Quot$ 
scheme parametrizing
subsheaves $\Cal K\subset \Cal F$ with $\frak H^{\Cal K} = \frak H$.  
We have then the bigraded versions of the ordinary and of the 
$C$-Grassmannian:
$$G_C(\frak H, P_{[(p_1,p_2),(q_1,q_2)]})\subset G(\frak H, 
P_{[(p_1,p_2),(q_1,q_2)]}),$$
and the following bigraded version of Theorem 1.4.1:

\proclaim{(4.4.3) Proposition} For $0\ll p_1+p_2\ll q_1+q_2$, the image of
the Grothendieck embedding $Sub_{\frak H}(\Cal F)\hookrightarrow
 G(\frak H, P_{[(p_1,p_2),(q_1,q_2)]})$ is 
$G_C(\frak H, P_{[(p_1,p_2),(q_1,q_2)]})$.
\endproclaim
The derived $C$-Grassmannian $RG_C(\frak H, P_{[(p_1,p_2),(q_1,q_2)]})$
is  defined as in
\S 3. It is a dg-manifold (up to quasiisomorphism), and we again set
$$RSub_{\frak H}(\Cal F)=RG_C(\frak H, P_{[(p_1,p_2),(q_1,q_2)]}).$$

Let now $\Cal K$ be a subsheaf of $\Cal F$ with  $\frak H^{\Cal K} = \frak H$.
 and let $Z$ be the connected component of $Sub(\Cal F)$ 
 containing $[\Cal K]$.
Set $h(s)=\frak H (s,0)$ and $k(t)= \frak H(0,t)$. Let $Y_1$ (respectively 
$Y_2$) be the
component of $RSub_h(\Cal F)$ (respectively $RSub_k(\Cal F)$) containing
$Z$ in the decomposition given by Proposition (4.4.2).

\proclaim{(4.4.4) Proposition} $Y_1$ and $Y_2$ are quasiisomorphic.
\endproclaim
\demo{Proof}
Via the obvious embeddings $A\hookrightarrow C$ and $B\hookrightarrow C$,
any $C$-module has also $A$- and $B$-module structures. This gives maps
$$RG_C(\frak H, P_{[(p_1,0),(q_1,0)]})\buildrel\varrho_1\over\longrightarrow
RG_A(h, M_{[p_1,q_1]}),$$
$$RG_C(\frak H, P_{[(0,p_2),(0,q_2)]})\buildrel\varrho_2\over\longrightarrow
RG_B(k, N_{[p_2,q_2]}).$$
By an obvious bigraded version of Theorem (4.3.2) $(a)$, 
$RG_C(\frak H, P_{[(p_1,0),(q_1,0)]})$ and 
$RG_C(\frak H, P_{[(0,p_2),(0,q_2)]})$
are quasiisomorphic.

In general, the maps $\varrho_1$ and $\varrho_2$ 
induce at the level 
of $\pi_0$ maps which ar 1--1, but not surjections, since the 
decomposition of $Sub(\Cal F)$ indexed by polynomials in two 
variables is finer than the one indexed by polynomials in one variable.  
Using (4.4.2), we can replace
$RG_C(\frak H, P_{[(p_1,0),(q_1,0)]})$ and 
$RG_C(\frak H, P_{[(0,p_2),(0,q_2)]})$
by quasiisomorphic  dg-schemes in which the connected components
of $\pi_0$ are ``separated''. Let $Z_1$ and $Z_2$ be the (quasiisomorphic)
respective components 
that contain $[\Cal K]$. By shrinking $Z_1$ and $Z_2$ if necessary, we get
induced maps $\varrho_i:Z_i\longrightarrow Y_i$, $i=1,2$.
The ``Whitehead theorem" (2.5.9) 
implies now that $\varrho_i$ are 
quasiisomorphisms. This completes the proof of Theorem 4.4.1.
\enddemo


\Refs
\widestnumber\key {HMS}

\ref\key BG \by A.K. Bousfield, V.K.A.M. Gugenheim \paper
On PL de Rham theory and rational homotopy type \jour Memoirs AMS
\vol 179 \yr 1976 \endref

\ref\key CK \by I. Ciocan-Fontanine, M. Kapranov \paper
Derived Hilbert schemes \jour in preparation \endref

\ref\key Fu\by W. Fulton \book Intersection Theory\publ
Springer-Verlag\yr 1984\endref

\ref\key Go \by G. Gotzmann \paper Eine Bedingung f\"ur die Flachheit
und das Hilbertpolynom eines graduierten Ringes \jour
Math. Zeit. \vol 158 \yr 1978\pages 61-70 \endref 


\ref \key Gr \by A. Grothendieck \paper Techniques \ de \  construction \ et 
\ th\'eor\`emes \ d'existence \ en \ g\'eometrie alg\'ebrique IV: Les sch\'emas
de Hilbert \yr 1960/61 \jour S\'eminaire Bourbaki \vol 221  \endref

\ref\key GV \by A. Grothendieck, J.-L. Verdier\paper Pr\'efaisceaux
(SGA4, Exp. I)  {\rm Lecture Notes in Mathematics 269}, \publ
Springer-Verlag\publaddr Berlin\yr 1972\endref


\ref\key Hi \by V. Hinich \paper Dg-coalgebras as formal stacks\jour
preprint math.AG/9812034 \endref

\ref\key HMS\by D. Husemoller, J. Moore, and J Stasheff
\paper Differential homological algebra and homogeneous spaces
\jour J. Pure Appl. Algebra \vol 5\yr 1975\pages 113-185\endref

\ref\key Ka1 \by M. Kapranov\paper Injective resolutions of BG and derived
moduli spaces of local systems\jour preprint alg-geom/9710027,
to appear in J. Pure Appl. Alg \endref

\ref\key Ka2 \by M. Kapranov \paper Rozansky-Witten invariants via Atiyah 
classes
\jour Compositio Math. \vol 115 \yr 1999\pages 71-113\endref

\ref\key Kol\by J. Koll\'ar \book Rational Curves on Algebraic Varieties
\publ Springer-Verlag \yr 1996\endref

\ref\key Kon \by M. Kontsevich \paper Enumeration of rational curves via
torus actions \jour in ``{\it The Moduli Space of Curves}" (R. Dijkgraaf, 
C. Faber, G. 
van der Geer Eds.) Progress in Math.\vol 129
\publ Birkhauser\publaddr Boston \yr1995 \pages 335-368 \endref

\ref\key Le \by D. Lehmann \paper Th\'eorie homotopique des
formes diff\'erentielles (d'apr\`es D. Sullivan) \jour Asterique 45, 
\yr 1977\endref

\ref\key Lo \by J.-L. Loday \book Cyclic Homology
\publ Springer-Verlag \yr 1995 \endref

\ref\key Ma \by M. Markl \paper A cohomology theory for A(m)-algebras and 
applications
\jour J. Pure Appl. Alg. \vol 83\yr 1992\pages 141-175\endref

\ref \key Mc\by J. McCleary\book User's Guide to Spectral Sequences,
{\rm Math. Lecture Ser., 12}\publ Publish or Perish Inc.
\publaddr Wilmington DE \yr 1985\endref

\ref \key Mu \by D. Mumford \book Lectures on curves on an algebraic
surface \publ Princeton Univ. Press \publaddr Princeton NJ \yr 1966 \endref

\ref\key Q1 \by D. Quillen \book Homotopical Algebra, {\rm Lecture Notes in
Math. 43}  \publ
Springer-Verlag\publaddr Berlin\yr 1967\endref

\ref\key Q2 \by D. Quillen \paper Rational Homotopy theory
\jour Ann. Math. \vol 90 \yr 1969 \pages 205-295\endref

\ref\key Re\by C. Rezk\paper Spaces of algebra structures and cohomology of 
operads
\jour thesis, MIT\yr 1996\endref 

\ref\key Se\by J.-P. Serre\paper Faisceaux alg\'ebriques
coh\'erents\jour Ann. of Math.\vol 61\yr 1955\pages 197-278
\endref

\ref\key Si \by C. Simpson \paper  Descente pour les $n$-champs\jour
preprint math.AG/9807049\endref

\ref\key St1\by J.D. Stasheff\paper Homotopy associativity of
H-spaces I, II\jour Trans. Amer. Math. Soc.\vol 108\yr1963\pages
275-312\endref

\ref\key St2 \by  J.D. Stasheff\paper Differential graded Lie algebras,
quasi-Hopf algebras and higher homotopy algebras 
\jour Quantum Groups (P.P. Kulish Ed.) (Lecture Notes Math. {\bf 1510}),
p. 120-137, Springer-Verlag, 1992 \endref

\ref\key Vi\by E. Viehweg \book Quasi-projective Moduli for Polarized
Manifolds,
{\rm Ergebnisse der Mathematik und ihrer Grenzgebiete (3), 30}\publ
Springer-Verlag\publaddr Berlin\yr 1995\endref

\endRefs
\enddocument